\documentclass[12pt,twoside,leqno]{article}
\usepackage{amsmath}
\usepackage{amssymb}
\usepackage{amsxtra}
\usepackage{amscd}
\usepackage{amsthm}
\usepackage[mathscr]{eucal}
\usepackage{color}
\usepackage{comment}

\theoremstyle{plain}
\newtheorem{thm}[subsection]{Theorem}
\newtheorem{prop}[subsection]{Proposition}
\newtheorem{cor}[subsection]{Corollary}
\newtheorem{lem}[subsection]{Lemma}

\theoremstyle{definition}
\newtheorem{defn}[subsection]{Definition}
\newtheorem{eg}[subsection]{Example}
\newtheorem{rem}[subsection]{Remark}
\newtheorem{para}[subsection]{}

\newtheorem{sbpara}[subsubsection]{}

\newenvironment{pf}{\proof[\proofname]}{\endproof}

\begin{document}

\title{Logarithmic geometry and Frobenius, II}

\author{Kazuya Kato, Chikara Nakayama, Sampei Usui}

\maketitle

\newcommand\Cal{\mathcal}
\newcommand\define{\newcommand}

\define\gp{\mathrm{gp}}%
\define\fs{\mathrm{fs}}%
\define\an{\mathrm{an}}%
\define\mult{\mathrm{mult}}%
\define\Ker{\mathrm{Ker}\,}%
\define\Coker{\mathrm{Coker}\,}%
\define\Hom{\mathrm{Hom}\,}%
\define\Ext{\mathrm{Ext}\,}%
\define\rank{\mathrm{rank}\,}%
\define\gr{\mathrm{gr}}%
\define\cHom{\Cal{Hom}}
\define\cExt{\Cal Ext\,}%

\define\cB{\Cal B}
\define\cC{\Cal C}
\define\cD{\Cal D}
\define\cO{\Cal O}
\define\cS{\Cal S}
\define\cM{\Cal M}
\define\cG{\Cal G}
\define\cH{\Cal H}
\define\cE{\Cal E}
\define\cF{\Cal F}
\define\cN{\Cal N}
\define\cY{\Cal Y}
\define\cT{\Cal T}
\define\fF{\frak F}
\define{\sW}{\Cal W}

\define\Dc{\check{D}}
\define\Ec{\check{E}}

\newcommand{\N}{{\mathbb{N}}}
\newcommand{\Q}{{\mathbb{Q}}}
\newcommand{\Z}{{\mathbb{Z}}}
\newcommand{\R}{{\mathbb{R}}}
\newcommand{\C}{{\mathbb{C}}}
\newcommand{\bN}{{\mathbb{N}}}
\newcommand{\bQ}{{\mathbb{Q}}}
\newcommand{\bF}{{\mathbb{F}}}
\newcommand{\bZ}{{\mathbb{Z}}}
\newcommand{\bP}{{\mathbb{P}}}
\newcommand{\bR}{{\mathbb{R}}}
\newcommand{\bC}{{\mathbb{C}}}
\newcommand{\bG}{{\mathbb{G}}}
\newcommand{\bbQ}{{\bar \mathbb{Q}}}
\newcommand{\ol}[1]{\overline{#1}}
\newcommand{\too}{\longrightarrow}
\newcommand{\respect}{\rightsquigarrow}
\newcommand{\compatible}{\leftrightsquigarrow}
\newcommand{\upc}[1]{\overset {\lower 0.3ex \hbox{${\;}_{\circ}$}}{#1}}
\newcommand{\Gmlog}{\bG_{m, \log}}
\newcommand{\Gm}{\bG_m}
\newcommand{\ep}{\varepsilon}
\newcommand{\Spec}{\operatorname{Spec}}
\newcommand{\nilp}{\operatorname{nilp}}
\newcommand{\prim}{\operatorname{prim}}
\newcommand{\val}{{\mathrm{val}}} 
\newcommand{\n}{\operatorname{naive}}
\newcommand{\bs}{\operatorname{\backslash}}
\newcommand{\Gal}{\operatorname{{Gal}}}
\newcommand{\gal}{{\rm {Gal}}({\bar \Q}/{\Q})}
\newcommand{\galp}{{\rm {Gal}}({\bar \Q}_p/{\Q}_p)}
\newcommand{\gall}{{\rm{Gal}}({\bar \Q}_\ell/\Q_\ell)}
\newcommand{\wep}{W({\bar \Q}_p/\Q_p)}
\newcommand{\wel}{W({\bar \Q}_\ell/\Q_\ell)}
\newcommand{\Ad}{{\rm{Ad}}}
\newcommand{\BS}{{\rm {BS}}}
\newcommand{\even}{\operatorname{even}}
\newcommand{\End}{{\rm {End}}}
\newcommand{\odd}{\operatorname{odd}}
\newcommand{\GL}{\operatorname{GL}}
\newcommand{\np}{\text{non-$p$}}
\newcommand{\g}{{\gamma}}
\newcommand{\G}{{\Gamma}}
\newcommand{\Lam}{{\Lambda}}
\newcommand{\La}{{\Lambda}}
\newcommand{\lam}{{\lambda}}
\newcommand{\la}{{\lambda}}
\newcommand{\uL}{{{\hat {L}}^{\rm {ur}}}}
\newcommand{\uQp}{{{\hat \Q}_p}^{\text{ur}}}
\newcommand{\sel}{\operatorname{Sel}}
\newcommand{\dt}{{\rm{Det}}}
\newcommand{\Sig}{\Sigma}
\newcommand{\fil}{{\rm{fil}}}
\newcommand{\SL}{{\rm{SL}}}
\renewcommand{\sl}{{\frak{sl}}}%
\newcommand{\spl}{{\rm{spl}}}
\newcommand{\st}{{\rm{st}}}
\newcommand{\Isom}{{\rm {Isom}}}
\newcommand{\Mor}{{\rm {Mor}}}
\newcommand{\bg}{\bar{g}}
\newcommand{\id}{{\rm {id}}}
\newcommand{\cone}{{\rm {cone}}}
\newcommand{\al}{a}
\newcommand{\ChL}{{\cal{C}}(\La)}
\newcommand{\Image}{{\rm {Image}}}
\newcommand{\toric}{{\operatorname{toric}}}
\newcommand{\torus}{{\operatorname{torus}}}
\newcommand{\Aut}{{\rm {Aut}}}
\newcommand{\Qp}{{\mathbb{Q}}_p}
\newcommand{\barQp}{{\mathbb{Q}}_p}
\newcommand{\Qpur}{{\mathbb{Q}}_p^{\rm {ur}}}
\newcommand{\Zp}{{\mathbb{Z}}_p}
\newcommand{\Zl}{{\mathbb{Z}}_l}
\newcommand{\Ql}{{\mathbb{Q}}_l}
\newcommand{\Qlur}{{\mathbb{Q}}_l^{\rm {ur}}}
\newcommand{\F}{{\mathbb{F}}}
\newcommand{\eps}{{\epsilon}}
\newcommand{\epsLa}{{\epsilon}_{\La}}
\newcommand{\epsLaVxi}{{\epsilon}_{\La}(V, \xi)}
\newcommand{\epsOLaVxi}{{\epsilon}_{0,\La}(V, \xi)}
\newcommand{\Qplin}{{\mathbb{Q}}_p(\mu_{l^{\infty}})}
\newcommand{\otimesQplin}{\otimes_{\Qp}{\mathbb{Q}}_p(\mu_{l^{\infty}})}
\newcommand{\galFl}{{\rm{Gal}}({\bar {\Bbb F}}_\ell/{\Bbb F}_\ell)}
\newcommand{\gallur}{{\rm{Gal}}({\bar \Q}_\ell/\Q_\ell^{\rm {ur}})}
\newcommand{\galFF}{{\rm {Gal}}(F_{\infty}/F)}
\newcommand{\galFv}{{\rm {Gal}}(\bar{F}_v/F_v)}
\newcommand{\galF}{{\rm {Gal}}(\bar{F}/F)}
\newcommand{\epsVxi}{{\epsilon}(V, \xi)}
\newcommand{\epsOVxi}{{\epsilon}_0(V, \xi)}
\newcommand{\plim}{\lim_
{\scriptstyle 
\longleftarrow \atop \scriptstyle n}}
\newcommand{\sig}{{\sigma}}
\newcommand{\ga}{{\gamma}}
\newcommand{\del}{{\delta}}
\newcommand{\Vss}{V^{\rm {ss}}}
\newcommand{\Bst}{B_{\rm {st}}}
\newcommand{\Dpst}{D_{\rm {pst}}}
\newcommand{\Dcrys}{D_{\rm {crys}}}
\newcommand{\DdR}{D_{\rm {dR}}}
\newcommand{\Fin}{F_{\infty}}
\newcommand{\Kla}{K_{\lambda}}
\newcommand{\Ola}{O_{\lambda}}
\newcommand{\Mla}{M_{\lambda}}
\newcommand{\Det}{{\rm{Det}}}
\newcommand{\Sym}{{\rm{Sym}}}
\newcommand{\LaSa}{{\La_{S^*}}}
\newcommand{\cX}{{\cal {X}}}
\newcommand{\MHG}{{\frak {M}}_H(G)}
\newcommand{\tauMla}{\tau(M_{\lambda})}
\newcommand{\Fvur}{{F_v^{\rm {ur}}}}
\newcommand{\Lie}{{\rm {Lie}}}
\newcommand{\cL}{{\cal {L}}}
\newcommand{\cW}{{\cal {W}}}
\newcommand{\fq}{{\frak {q}}}
\newcommand{\cont}{{\rm {cont}}}
\newcommand{\SC}{{SC}}
\newcommand{\Om}{{\Omega}}
\newcommand{\dR}{{\rm {dR}}}
\newcommand{\crys}{{\rm {crys}}}
\newcommand{\hatSig}{{\hat{\Sigma}}}
\newcommand{\rdet}{{{\rm {det}}}}
\newcommand{\ord}{{{\rm {ord}}}}
\newcommand{\BdR}{{B_{\rm {dR}}}}
\newcommand{\BdRO}{{B^0_{\rm {dR}}}}
\newcommand{\Bcrys}{{B_{\rm {crys}}}}
\newcommand{\Qw}{{\mathbb{Q}}_w}
\newcommand{\barkappa}{{\bar{\kappa}}}
\newcommand{\cP}{{\Cal {P}}}
\newcommand{\cZ}{{\Cal {Z}}}
\newcommand{\cA}{{\Cal {A}}}
\newcommand{\oppLa}{{\Lambda^{\circ}}}

\renewcommand{\bar}{\overline}
\newcommand{\et}{\mathrm{\acute{e}t}}
\newcommand{\loget}{\mathrm{log\acute{e}t}}
\newcommand{\pri}{{{\rm prim}}}
\newcommand{\add}{{\rm{add}}}

\define\cR{\Cal R}
\newcommand{\br}{{\bold r}}
\newcommand{\mild}{{{\rm {mild}}}}
\newcommand{\nar}{{{\rm {nar}}}}


\begin{abstract} 
  Based on the strong analogy between the category of log mixed Hodge structures and the category $\cA_X$ of $\ell$-adic nature, which we have introduced in the previous part and is closely related to the weight-monodromy conjecture, we prove the $\ell$-adic analogues of some theorems in Hodge theory related to the $\SL(2)$-orbit theorem. 
\end{abstract}

\renewcommand{\thefootnote}{\fnsymbol{footnote}}


\section*{Contents}

\noindent 
Introduction

\noindent \S\ref{s:ms}. Monodromy systems 

\noindent \S\ref{s:ratio}. Ratio spaces

\noindent \S\ref{s:na}. Analogy

\noindent \S\ref{s:twi}. SL(2)-orbit theorems on the torus twist back

\noindent \S\ref{s:1var}. SL(2)-orbit theorems in one variable

\noindent \S\ref{s:spl}. SL(2)-orbit theorems on $\spl_W$

\noindent \S\ref{s:del}. SL(2)-orbit theorems on $\delta_W$

\noindent \S\ref{s:reg}. Asymptotic behaviors of regulators and local height pairings

\noindent \S\ref{s:linalg}. On the expectation 6.16 of Part I

\noindent References

\section*{Introduction}\label{s:intro}

\begin{para}
In our previous paper \cite{KNU} (which we call Part I), 
we have introduced an abelian category $\cA_X$ of $\ell$-adic sheaves with weight filtrations on an fs log scheme $X$ over a finite field.
There is a strong analogy 
between this category and the category of log mixed Hodge structures on an fs log analytic space $X$ over $\C$: 

\smallskip

$\qquad$ $\cA_X$ on $X$ over a finite field $\quad \longleftrightarrow\quad$  LMH on $X$ over $\C$ 

\smallskip

\noindent Our guiding principle is that every theorem on the one side should always have an analogue on the other side. 
  It gives a new relation between $\ell$-adic theory and Hodge theory. 
  Based on this principle, in this paper, we will obtain various $\SL(2)$-orbit theorems on the left-hand-side which are analogous to the $\SL(2)$-orbit theorems on the right-hand-side.

  At each point of $X$, by taking fibers of both sides of the above analogy, 
another analogy 
$$\text{monodromy system}\;\;\longleftrightarrow\;\;\text{nilpotent orbit}$$
appears. SL(2)-orbit theorems in Hodge theory are theorems about nilpotent orbits. The SL(2)-orbit theorems in this paper are theorems about monodromy systems. 
Roughly speaking, a nilpotent orbit is $(N_1, \dots, N_n, F)$, where $N_1, \dots, N_n$ are monodromy operators and $F$ is a Hodge filtration satisfying certain conditions. Roughly speaking, a monodromy system is a family of monodromy operators (see Definition \ref{na0} for the precise definition of monodromy system) satisfying certain conditions. As is explained in \ref{analogy}, we have an analogy
$$\text{monodromy system $(N_1, \dots, N_{n+1})$}\;\;\longleftrightarrow\;\;\text{nilpotent orbit $(N_1, \dots, N_n, F)$}.$$
One of the $\SL(2)$-orbit theorems obtained  in this paper is the following Theorem \ref{thm2}, which is an analogue of the following Theorem \ref{thm1}  in Hodge theory proved in \cite{KNU08}. (These Theorems \ref{thm1} and \ref{thm2} are stated here in slightly rough forms.)

\end{para}
  
  \begin{thm}\label{thm1} (Cf.\ {\rm \cite{KNU08}} Theorem $0.5$.) Let $(N_1, \dots, N_n,F)$ be a nilpotent orbit. Then for $y_j>0$ ($1\leq j\leq n$), $y_j/y_{j+1}\gg 0$ ($y_{n+1}$ denotes $1$), $\spl_W(\exp(\sum_{j=1}^n iy_jN_j)F)$ is expressed by  convergent Taylor series in $y_{j+1}/y_j$ ($1\leq j\leq n$).

  \end{thm}
  
  \begin{thm}\label{thm2} (Cf.\ Theorem $\ref{splthm}$ in this paper.) Let $(N_1, \dots, N_{n+1})$ be a monodromy system. Then for $y_j>0$ ($1\leq j\leq n$), $y_j/y_{j+1}\gg 0$ ($y_{n+1}$ denotes $1$), $\spl_W((\sum_{j=1}^n iy_jN_j)+N_{n+1})$ is expressed by  convergent Taylor series in $y_{j+1}/y_j$ ($1\leq j\leq n$).

  \end{thm}
  
  Here in these theorems, $\spl_W(\bullet)$ is the canonical splitting of the weight filtration $W$ associated to $\bullet$.
  
  \begin{para}
  The asymptotic behaviors of the invariants  $\spl_W$ and $\delta_W$  for a nilpotent orbit are the  main subjects in \cite{KNU08}, \cite{KNU4}, etc. The asymptotic behaviors of the invariants  $\spl_W$ and $\delta_W$  for a monodromy system are the  main subjects of this paper.
  We apply such results to study asymptotic behaviors 
of regulator maps and local height pairings, 
which appear in the algebraic geometry over non-archimedean local fields (see Section \ref{s:reg}). 
  Further, we give some complementary results to Part I (see Section \ref{s:linalg}).

  The organization of this paper is as follows. 
  In Sections \ref{s:ms} and \ref{s:ratio}, we introduce monodromy systems 
and 
ratio spaces, respectively. 
  In Section \ref{s:na}, we explain the above analogy more in detail. 
  In Section \ref{s:twi}, we prove refinements of Theorem 5.3 of Part I which are used in Sections \ref{s:1var}--\ref{s:del}. 
  In Section \ref{s:1var}, we prove the one-variable cases of the main results to which the general cases are reduced. 
  In Sections \ref{s:spl} and \ref{s:del}, we prove two main theorems in this paper (Theorems \ref{splthm} and \ref{delthm0}) concerning the asymptotic behaviors of  $\spl_W$ and  $\delta_W$ of monodromy systems, respectively. 
  Section \ref{s:reg} is an application to geometry over a non-archimedean local field. 
  In Section \ref{s:linalg}, we give a simple example on which Question 6.3 of Part I is affirmative. (This section is independent of the previous ones.)   
  \end{para}
  
  \begin{para}

We hope to remark again that as was explained in Remark 1.13 of Part I, one of the authors (K. Kato) published a wrong paper \cite{Ka:KJM} on the subject SL(2)-orbit theorems related to Deligne systems, and  hopes that  \cite{Ka:KJM} is never used by any person. The present paper is a correct work on the same subject.

\end{para}

\begin{para}

The authors thank Spencer Bloch for his stimulation to this work. Some parts of this paper were originally prepared for the joint works of the first author  with him on the asymptotic behaviors of regulators and height pairings. 

\end{para}
  K.\ Kato was 
partially supported by NFS grants DMS 1601861 and DMS 2001182.
C.\ Nakayama was 
partially supported by JSPS Grants-in-Aid for Scientific Research (C) 16K05093 and (C) 21K03199.
S.\ Usui was 
partially supported by JSPS Grants-in-Aid for Scientific Research (C) 17K05200. 

\section{Monodromy systems}\label{s:ms}
We prove some basic facts on monodromy systems and associated Deligne systems.

\begin{para}\label{fgc}
 By a {\it finitely generated cone}, we mean a commutative monoid $\sig$ (whose semi-group law is written  additively) endowed with an action $\R_{\geq 0}\times \sig \to \sig$ of the multiplicative monoid $\R_{\geq 0}= \R^{\mult}_{\geq 0}=\{r\in \R\;|\; r\geq 0\}$ such that as a monoid with an action of $\R_{\geq 0}$, 
$\sig$ is isomorphic to $\R_{\geq 0} v_1+\dots +\R_{\geq 0} v_n$ for some $\R$-vector space $V$ and for some 
 $v_1, \dots, v_n\in V$.

  For a finitely generated cone $\sig$, a {\it face} of $\sig$ is a submonoid $\tau$ of $\sig$ such that if $a, b\in \sig$ and $a+b\in \tau$, then $a,b\in \tau$. 
  Then any face is stable under the action of $\R_{\geq 0}$ and it itself is a finitely generated cone. 
Let $\text{face}(\sig)$ be the set of all faces of $\sig$. It is a finite set. For a homomorphism $\sig\to \sig'$ of finitely generated cones, we have the map $\text{face}(\sig)\to \text{face}(\sig')$ which sends $\tau\in \text{face}(\sig)$ to the smallest face of $\sig'$ which contains the image of $\tau$.
 We also have a map $\text{face}(\sig')\to \text{face}(\sig)$ which sends $\tau' \in$ face$(\sig')$ to its inverse image. 
 
  The group $\sig^{\gp}=\{x-y\,|\, x,y\in \sig\}$ associated to a finitely generated cone $\sig$ is 
 also written as $\sig_\R$. It has a canonical structure of a finite dimensional  $\R$-vector space. 
   
 A finitely generated cone is said to be {\it sharp} if $\sig \cap (-\sig)=\{0\}$ in $\sig_\R$. 
 
For a finitely generated  cone $\sig$, the set $\sig\spcheck$ of homomorphisms from $\sig$ to $\R^{\add}_{\geq 0}=\R_{\geq 0}$ (regarded as a finitely generated  cone) is  a sharp finitely generated cone. We have $(\sig\spcheck)\spcheck=\sig$ if $\sig$ is sharp. 
\end{para}

\begin{defn}\label{na0} A {\it monodromy system} is a quadruple $(V, W, \sig, Y)$, where 
  $V$ is  a finite dimensional 
$\C$-vector space, $W$ is a finite increasing filtration on $V$, $\sig$ is a sharp finitely generated cone with an admissible action on $V$  (\cite{KNU} 2.3), and $Y$
 is a splitting of $W(\sig)$ which is compatible with $W$ and for which all elements of $\sig$ are of weight $-2$. 
 
 Here the splitting $Y$ is formulated as an endomorphism of $V$ whose all eigenvalues are integers such that $W(\sig)_w= \bigoplus_{w'\leq w}\; \text{Ker}(Y-w')$ for all $w$.

  If $(V,W,\bR^n_{\ge0}, Y)$ is a monodromy system, where the standard generators of $\bR^n_{\ge0}$ act via $N_j:V\to V$ ($1\leq j\leq n$), we say that 
$(V,W,N_1,\dots,N_n, Y)$ is a monodromy system. 
 Here $V, W$, and $Y$ can be omitted, and 
we often simply say that $(N_1, \dots, N_n)$ is a monodromy system. 

\end{defn}

\begin{para}

Let $X$ be $\Spec(k)$ for a finite field $k$ endowed with an fs log structure, and let $H$ be an object of $\cA_X$ (\cite{KNU} 1.1). Fix a lifting of Frobenius $\in \Gal(\bar k/k)$ to $\pi_1^{\log}(X)$. Then the stalk of  $H$ is regarded as a monodromy system. 

This is similar to the fact that if $X$ is an fs log point over $\C$, an LMH on $X$ is regarded as a nilpotent orbit (\cite{KNU4} 2.2.2). 

\end{para}

\begin{para} The category of monodromy systems for a fixed $\sig$ is an abelian category. This is shown by the same arguments in \cite{KNU} Section 4 in which we proved that  $\cA_X$ is an abelian category for an fs log scheme $X$ of finite type over a finite field.

\end{para}

\begin{lem}\label{tauMS} Let $(V, W, \sig, Y)$ be a monodromy system. Then for a face $\tau$ of $\sig$, $(V, W(\tau), \sig, Y)$ is a monodromy system.

\end{lem}

\begin{pf} The admissibility of the action of $\sig$ on $(V, W(\tau))$ follows from the admissibility of the action of $\sig$ on $(V, W)$ (for a face $\tau'$ of $\sig$, $W(\tau')$ for the former is $W(\tau'')$ for the latter, where $\tau''$ is the smallest face of $\sig$ which contains $\tau\cup \tau'$).  For any element $N$ of the interior of $\tau$, 
 $W(\tau)$ is the relative monodromy filtration of $N$ with respect to $W$ 
 and  $N$ is of weight $-2$ for $Y$. Hence $Y$ is compatible with $W(\tau)$. 
\end{pf}

\begin{prop}\label{adm} 
Let $\sig\to \sig'$ be a homomorphism of sharp finitely generated  cones,  and assume that $\sig'$ acts on a finite dimensional $\C$-vector space $V$ endowed with a finite increasing filtration $W$. 

$(1)$ If the action of $\sig'$ on $V$ is admissible, the induced action of $\sig$ on $V$ is admissible. For a face $\tau$ of $\sig$, $W(\tau)$ coincides with $W(\tau')$, where $\tau'$ is the image of $\tau$ under the map $\mathrm{face}\,(\sig)\to \mathrm{face}\,(\sig')$ ($\ref{fgc}$). 

$(2)$ If the induced action of $\sig$ on $V$ is admissible and if the map $\sig\to \sig'$ is surjective, the action of $\sig'$ on $V$ is admissible. 
For a face $\tau'$ of $\sig'$, $W(\tau')$ coincides with $W(\tau)$, where $\tau$ is the image of $\tau'$ under the map $\mathrm{face}\,(\sig')\to \mathrm{face}\,(\sig)$ ($\ref{fgc}$). 
\end{prop}

\begin{pf}
  The proof is similar to that for \cite{KNU3} Lemma 1.2.3, though 
the definition of the admissibility there is slightly different from the one in this paper. 
\end{pf}

\begin{cor}
\label{c:adm}
$(1)$ Let $(V, W, \sig,Y)$ be a monodromy system. If $N_1, \dots, N_n$ are  elements of $\sig$ such that $\sig$ is the smallest face of $\sig$ which contains them,  then  $(N_1, \dots, N_n)$ (with $V$, $W$, the cone $\R^n_{\geq 0}$ acting on $V$ via $(N_j)_{1\leq j\leq n}$, and $Y$) is a monodromy system.

$(2)$ Let $(N_1, \dots, N_n)$ be a monodromy system (with the cone $\R^n_{\geq 0}$ acting on $V$ via $(N_j)_j$) and let $\sig$ be the image of the homomorphism $\R^n_{\geq 0}\to \End(V)$. Assume that $\sig$ is sharp. Then  $(V, W,\sig, Y)$ is a monodromy system. 
\end{cor}

\begin{pf} In fact, (1) follows from (1) of Proposition \ref{adm} and (2) follows from (2) of Proposition \ref{adm}. \end{pf}

\begin{para}\label{DS} We consider how Deligne systems  (\cite{SC} Section 2,  \cite{BPR} Section 6,  reviewed in Part I \cite{KNU} Section 2) are compared with and related to monodromy systems.
In this section, we consider Deligne systems over $\C$ though Deligne systems can be considered over any coefficient field.  
We will often call a Deligne system $(V, (W^j)_{0\leq j\leq n}, (N_j)_{1\leq j\leq n}, Y)$ in $n$ variables  simply a Deligne system $(N_1, \dots, N_n)$ or a Deligne system $(N_1, \dots, N_n, Y)$.

Recall that for a Deligne system $(N_1, \ldots, N_n,Y)$ and 
 for $0\le j \le n$, we have the splitting $Y^j$ of $W^j$ by the downward induction:  $Y^n=Y$ and 
$Y^{j-1}$ is the Deligne  splitting (\cite{KNU} 2.7) of $W^{j-1}$ obtained by $N_j$ and  $Y^j$.
\end{para}

\begin{prop}
\label{p:ms&Ds}
  Let $(V, W, \sig, Y)$ be a monodromy system and let $N_1, \dots, N_n$ be elements of $\sig$ such that $\sig$ is the smallest face of $\sig$ which contains them. 
  Then we have a Deligne system $(V, (W^j)_j, (N_j)_j, Y)$, where $W^j=W(\sig_j)$ for $0\le j \le n$ with $\sig_j$ being the smallest face of $\sig$ containing $N_1,\ldots,N_j$ (in particular, $W^0=W$ and $W^n=W(\sig)$).
  
\end{prop}

\begin{pf}
  Let $0\le i < j \le n$. 
  Let $U$ be a filter of $W(\sig_i)$.
  By admissibility, we have that 
$W(\sig_j)\vert_{U}=M((N_1+\cdots+N_j)|_{U}, W(\sig_{j-1})|_{U})=M(N_j|_{U}, W(\sig_{j-1})|_U)$. 
  In particular, 

\smallskip

\noindent $(*)$ $W(\sig_j)=M(N_j, W(\sig_{j-1}))$

\smallskip
\noindent  
for $0<j\le n$.  
  Since $Y$ is compatible with $W=W(\sig_0)$ and all $N_j$ are of pure weight for $Y$, 
$(*)$ inductively implies that $Y$ is compatible with the other $W(\sig_j)$.
  Finally, again by the admissibility, we have that $N_i$ preserves $W(\sig_j)$ for any $i$ and $j$, and that $N_i$ belongs to $W^j_{-2}\End(V)$ for $i \le j$. 
\end{pf}

\begin{lem}\label{subDS} Let $(N_1, \dots, N_n, Y)$ be a Deligne system. Then 
 for $1 \le j \le n$, 
$(N_1,\cdots, N_j, Y^j)$ is a Deligne system.  

\end{lem}

\begin{pf} This is  a part of \cite{SC} Proposition 2. \end{pf}

An analogue of Lemma \ref{subDS} for monodromy systems is the following 
Lemma \ref{subMS}.

\begin{lem}\label{subMS} 
Let $(V, W, \sig, Y)$ be a monodromy system, let $\tau$ be a face of $\sig$, and let $N$ be an element of $\sig$ whose image in $\sig/\tau:=(\sig +\tau_\R)/\tau_\R$  belongs to the interior of $\sig/\tau$. Then $\sig$ is the smallest face of $\sig$ which contain $\tau$ and $N$, and we have the Deligne system $(V, (W(\tau), W(\sig)), N, Y)$ in one variable. For the associated Deligne splitting $Y'$ of $W(\tau)$, $(V, W, \tau, Y')$ is a monodromy system.

\end{lem}

\begin{pf} The first part follows from Lemma \ref{tauMS} and Proposition \ref{p:ms&Ds}. The proof of the second part is similar to that of Lemma \ref{subDS}. The crucial point is that all $N'\in \tau$ are of weight $-2$ for $Y'$, which is proved as follows. We have the monodromy system $(V, W(\tau), \sig, Y)$ by Lemma \ref{tauMS}, and for the pair $((V, W(\tau), \sig, Y),N)$ of the monodromy system $(V, W(\tau), \sig,Y)$ and $N$, $N'$ is a homomorphism of pairs $((V, W(\tau), \sig, Y), N)\to ((V, W(\tau), \sig, Y), N)(-1)$, where $(-1)$ is the Tate twist. Hence $N'$ gives a homomorphism of Deligne systems $(V, (W(\tau), W(\sig)), N, Y)\to (W, (W(\tau), W(\sig)), Y)(-1)$.
By the functoriality of the construction of the Deligne splitting,  $N'$ induces a morphism of triples $(V, W(\tau), Y') \to (V, W(\tau), Y')(-1)$ and hence $N'$ is of weight $-2$ for $Y'$. \end{pf} 

\begin{para}\label{spldel} A Deligne system in one variable is equivalent to a monodromy system with $\sig=\R_{\geq 0}$.

For a Deligne system $(V, (W, W'), N, Y)$ in one variable, we denote by $\spl_W(N)$ the Deligne splitting $Y^0$ of $W$. We define $\delta_W(N)\in \gr^W\End(V)$ as the map $\gr^WV \to\gr^WV$ corresponding to $N-N^{[0]}:V \to V$, where  $N^{[0]}$ is the component of $N$ of $Y^0$-weight $0$, via the isomorphism $\gr^WV\cong V$ given by 
$Y^0$ of $W$.
 $N$ is recovered from $\gr^W(N)$, $\spl_W(N)$ and $\delta_W(N)$: $N$ is the map $V\to V$ corresponding to  $\gr^W(N)+\delta_W(N): \gr^WV\to \gr^WV$ via the splitting $Y^0$ of $W$.

Properties of these
$$\spl_W(N) \quad \text{and}\quad \delta_W(N)$$
are the main subjects of this paper.

Note that  $\delta_W(N)\in  \gr^W\End(V)$ is independent of the choice of $Y$ (\cite{Ka1} Proposition 1.5.7). That is, $\delta_W(N)$ is defined canonically for $(V, W, N)$ such that the relative monodromy filtration $M(W, N)$  exists and such that there is a splitting $Y$ of $M(W, N)$ which is compatible with $W$ and for which $N$ is of weight $-2$. 
  
 Thus if $(V, W, \sig, Y)$ is a monodromy system, $\delta_W(N)\in \gr^W\End(V)$ is defined canonically for all $N\in \sig$ because such a splitting of $M(W, N)$ exists for the following reason. Take an element $N'$ of the interior of $\sig$. Then we have a Deligne system $(V, (W^j)_{0\leq j\leq 2}, (N, N'), Y)$ ($W^0=W, W^1=M(W, N), W^2=W(\sig)$)  by 
Proposition \ref{p:ms&Ds} and the Deligne splitting $Y^1$ of  $M(W, N)$ is such a splitting (Lemma \ref{subMS}).

\end{para}

\begin{para}\label{sdwt} The following (1)--(4) will be used in later sections.

(1) If $W$ is a finite increasing filtration of a finite dimensional vector space $V$ and if $Y$ and $Y'$ are splittings of $W$, there exists a unique $u\in \Aut(V)$ such that $Y'=uYu^{-1}$ and $u-1\in W_{-1}\End(V)$. 

This is seen by the fact that  if $s, s':\gr^WV \overset{\cong}\to V$ are the isomorphisms given by $Y$ and $Y'$, respectively, $Y'=uYu^{-1}$ is equivalent to $s'=u\circ s$. 

In the following (2)--(4), let $(V, (W, W'), N, Y)$ be a Deligne system in one variable. 

(2) Let $N_*\in\End(V)$ and assume that $(V, (W, W'), N_*, Y)$ is also a Deligne system in one variable. Let $Y^0=\spl_W(N)$, $Y^0_*=\spl_W(N_*)$ and write $Y^0=uY^0_*u^{-1}$ 
with $u\in \Aut(V)$ such that $u-1\in W_{-1}\End(V)$. Then $u$ is of $Y$-weight $0$.

This follows from the fact that $Y^0$ and $Y^0_*$ are of $Y$-weight $0$ (that is, they are compatible with $Y$) and from the uniqueness of $u$.

(3) $\delta_W(N)$ belongs to $W_{-2}\gr^W\End(V)= \bigoplus_{w\leq -2} \gr^W_w\End(V)$.

This is seen as follows. Let $N^{[-d]}$ be the component of $N$ of $Y^0$-weight $-d$: $N=\sum_{d \geq 0}N^{[-d]}$.
  Then, by the definition of Deligne splitting, for every $d>0$, 
$N^{[-d]}$ is killed by $\Ad(N^{[0]})^{d-1}$. 
  In particular, $N^{[-1]}=0$.

(4) $\delta_W(N)$ is of $Y$-weight $-2$. 

In fact, if we identify $\gr^WV$ and $V$ by the Deligne splitting, we have $N=\gr^W(N)+\delta_W(N)$. Since $N$ and $\gr^W(N)$ are of $Y$-weight $-2$, so is $\delta_W(N)$.

\end{para}

\begin{para}\label{calD} Let $V$ be a finite dimensional vector space over $\C$ endowed with finite increasing filtrations $W$ and $W'$ and with a splitting $Y$ of $W'$ which is compatible with $W$. Let $\overline{\cD}$ be the set of all $N\in \End(V)$ which are of weight $-2$ for $Y$, and let $\cD$ be the subset of $\overline{\cD}$ consisting of all elements $N$ such that $W'$ is the relative monodromy filtration of $N$ with respect to $W$. Then $\overline{\cD}$ is a closed complex analytic submanifold of $\End(V)$ and $\cD$ is an open set of $\overline{\cD}$. We have complex analytic maps 
$$\spl_W: \cD \to \spl(W), \quad \delta_W: \cD \to \gr^W\End(V),$$
where $\spl(W)$ is the complex analytic manifold of all splittings of $W$. 

This set $\cD$ is an analogue of the classifying space of mixed Hodge structures (in the case where $V$ is $W$-pure, it is an analogue of the Griffiths period domain). The Hodge analogues of the maps $\spl_W$ and $\delta_W$ are discussed, for example, in \cite{KNU4} Section 1.2.

\end{para}

\begin{para}\label{SL2}
Let $(N_1, \dots, N_n)$ be a Deligne system. 

For $0\leq j\leq n$, let $\tau_j$ be the action of $\Gm$ corresponding to $Y^j$.

  For $1 \le j \le n$, let $\hat N_j$ be the component of $N_j$ of weight $0$ for the splitting $Y^{j-1}$ of $W^{j-1}$. Then $\hat N_j$ is of weight $0$ for the splitting $Y^i$ of $W^i$ for $0 \le i <j$
(\cite{SC} Corollary 1 to Proposition 2).

As in \cite{SC}, there is an action of $\mathbb{G}_m \times \SL(2)^n$ on $V$ characterized as follows. 

The action of the $\mathbb{G}_m$ is given by $\tau_0$. For $1\leq j\leq n$,  $$\tau_j\tau_{j-1}^{-1}$$ coincides with the composition $\mathbb{G}_m\to \SL(2)^n\to \Aut(V)$, where the $k$-th component 
($1\leq k\leq n$) of the first arrow is $a \mapsto \begin{pmatrix} a^{-1} & 0\\ 0 & a\end{pmatrix}$ for $k=j$ and is trivial otherwise. For $1\leq j\leq n$, $\begin{pmatrix} 0&1\\ 0&0\end{pmatrix}$ in the $j$-th $\sl(2)$ in $\sl(2)^n$ acts as  $\hat N_j$.

  For a monodromy system $(N_1, \ldots, N_n)$, we define $Y^j$, 
$\tau_j$  
$(0\le j \le n)$ and $\hat N_j$ $(1 \le j \le n)$ as those for the associated Deligne system (Proposition \ref{p:ms&Ds}).

\end{para}

 \begin{para}\label{SLrep} Recall that all finite dimensional representations of $\SL(2)^n$ are semi-simple and every irreducible representation of $\SL(2)^n$ has the shape
 $$\Sym^{m(1)}V_1 \otimes \dots\otimes  \Sym^{m(n)}V_n,$$ where $m(i)\geq 0$ and $V_i$ is the $2$-dimensional vector space on which $\SL(2)^n$ acts via its $i$-th component in the standard way. 
 
 \end{para}

\begin{prop}\label{MSD2} Let $(N_1, \dots, N_n)$ be a Deligne system. Then $(N_1, \hat N_2, \dots, \hat N_n)$ is a monodromy system, and $W(\R^n_{\geq 0})$ of it coincides with $W^n$ of the Deligne system.

\end{prop}

\begin{pf} This follows from \ref{SLrep}.  \end{pf}

\begin{prop}\label{jk1} Consider a Deligne system $(N_1, \dots, N_n)$ and the associated representation of $\mathbb{G}_m \times \SL(2)^n$ ($\ref{SL2}$). For each $a\in \Z$, let $N_j^{[a]}$ be the component of $Y^{j-1}$-weight  $a$ of  $N_j$. Then 
$[N_j^{[a]}, \hat N_k]=0$ if $1\leq j< k$. 
\end{prop}

\begin{pf}  Let $N_k^{(a)}$ be the component of $N_k$ of $\tau_{k-1}$-weight $a$. Since $N_j$ is of $\tau_{k-1}$-weight $-2$, applying $\Ad(\tau_{k-1}(c))$ to $[N_j, N_k]=0$, we have $c^{-2}\sum_a c^a[N_j, N_k^{(a)}]=0$. Hence $[N_j, N_k^{(a)}]=0$ for all $a$. By taking $a=0$, we have $[N_j, \hat N_k]=0$. The $\tau_{j-1}$-weight of $\hat N_k$ is $0$. Applying $\Ad(\tau_{j-1}(c))$ to  $[N_j, \hat N_k]=0$, we have 
$\sum_a c^a[N_j^{[a]}, \hat N_k]=0$. Hence $[N_j^{[a]}, \hat N_k]=0$ for all $a$. \end{pf}

\begin{prop}\label{twospl} Let $(N_1, \dots, N_n)$ be a Deligne system. Then the splitting $Y^0$ of $W$ coincides with 
$\spl_W(N_1+ \sum_{j=2}^n \hat N_j)$.  

\end{prop}

\begin{pf}  Let $(\,\cdot\,)^{[w]}$ be the component of weight $w$ for the splitting $Y^0$ of $W$. It is sufficient to prove that for every $w<0$, $(N_1+\sum_{j=2}^n \hat N_j)^{[w]}$ belongs to the primitive part for $\sum_{j=1}^n \hat N_j$ acting on $\gr^W_w$. Let $\mu=-w$. It is sufficient to prove that $\Ad(\sum_{j=1}^n \hat N_j)^{\mu-1}((N_1+\sum_{j=2}^n \hat N_j)^{[w]})=0$. We have $(N_1+\sum_{j=2}^n  \hat N_j)^{[w]}=N_1^{[w]}$. Furthermore, by Proposition \ref{jk1},  $[\hat N_j, N_1^{[w]}]=0$ for $2\leq j\leq n$. Hence we are reduced to $\Ad(\hat N_1)^{\mu-1}(N_1^{[w]})=0$. But the last thing follows from the fact that $Y^0$ is the Deligne splitting of $W$ for $W^1$ and $N_1$.
\end{pf}

\begin{prop}\label{jk2} 
  Let the situation be as in Proposition $\ref{jk1}$. 
  Let $1\leq j\leq k$ and $a \in \bZ$. 
  Let $\hat N_k^+$ be the action of $\begin{pmatrix} 0& 0\\ 1& 0\end{pmatrix}$ in the $k$-th $\sl(2)$ in $\sl(2)^n$. 
  Then we have $[N_j^{[a]}, \hat N_k^+]=0$ unless $j=k$ and $a=0$. 
\end{prop}

\begin{pf} Assume first $a\neq 0$. For $j\leq \ell\leq  k$, by the fact that $W^{\ell}$ is the relative monodromy filtration of $N:=N_j + \sum_{j<j'\leq \ell} \hat N_{j'}$ with respect to $W^{j-1}$ and by the fact that $Y^{j-1}$ gives the Deligne splitting of $W^{j-1}$ associated to the splitting $Y^{\ell}$ of $W^{\ell}$ with respect to $N$ by Proposition \ref{twospl}, we have $[N^{[a]}, \sum_{j\leq j'\leq \ell} \hat N_{j'}^+]=0$. Since $\hat N_{j'}^{[a]}=0$ for $j'>j$, we have $N^{[a]}= N_j^{[a]}$, and hence we have $[N_j^{[a]}, \sum_{j\leq j'\leq \ell} \hat N_{j'}^+]=0$. This proves $[N_j^{[a]}, \hat N_k^+]=0$ for $k\geq j$.   

Assume next $a=0$. Then $N_j^{[a]}=\hat N_j$ and $[\hat N_j, \hat N_k^+]=0$ if $j\neq k$. 
\end{pf}

\section{Ratio spaces}\label{s:ratio}

For a sharp finitely generated cone $\sig$, we define and discuss the space $\sig_{[:]}$ of ratios. This story of the space of ratios is essentially a rewriting of the story of the space of ratios in \cite{KNU4} Section 4 and \cite{KNU} Section 3. In fact, for a sharp fs monoid $\cS$ and the dual cone $\sig=\Hom(\cS, \R_{\geq 0}^{\add})$, where $\R_{\geq 0}^{\add}$ denotes the additive monoid $\R_{\geq 0}= \{x\in \R\;|\; x\geq 0\}$, $\sig_{[:]}$ here coincides with the space of ratios  $R(\cS)$ defined in \cite{KNU4} 4.1.3 and reviewed in \cite{KNU} 3.1. The structure of a real analytic manifold on $\sig_{[:]}$ discussed in this Section is essentially considered in \cite{KNU6}  1.4.19 (see \ref{Aan1}). 
Thus we are only presenting the story in \cite{KNU4}, \cite{KNU}, and \cite{KNU6} in a dual way. But this formulation has some advantages for the presentations of the results of this paper. 
 The proofs in this section are often omitted if they are essentially given in \cite{KNU4} or \cite{KNU6}.

\begin{para} 
For a sharp finitely generated  cone $\sig$, the set of ratios $\sig_{[:]}$ associated to $\sig$ is the set of maps $r: \sig\spcheck\times \sig\spcheck\smallsetminus \{(0,0)\}\to [0, \infty]$ satisfying the following conditions (i)--(iv).

(i) $r(f,g)=r(g,f)^{-1}$.

(ii) $r(f,g)r(g,h)= r(f,h)$ if $\{r(f,g), r(g,h)\}\neq \{0,\infty\}$.

(iii) $r(f+g, h)=r(f,h)+ r(g, h)$. 

(iv) $r(cf, g)= cr(f,g)$ for $c\in \R_{>0}$.

We endow $\sig_{[:]}$ with the topology of simple convergence.

\end{para}

\begin{eg}\label{eg}
Let $\sig=\R_{\geq 0}^2$. Then we have a homeomorphism  $$\sig_{[:]}\cong [0,\infty]\;\; ;\;\;  r\mapsto r(t_1,t_2),$$
where $t_j$ are the standard coordinate functions of $\sig$. 
\end{eg}

\begin{para} For an fs log point $x$, the space of ratios $x_{[:]}$ in \cite{KNU} Section 3 coincides with  $\sig_{[:]}$ for the monodromy cone $\sig=  \Hom((M_x/\cO^\times_x)_{\bar x}, \R_{\geq 0}^{\add})$ of $x$ ($\bar x$ denotes the separable closure of $x$). 
\end{para}

\begin{para}\label{1:1} 

 There is a bijection between $\sig_{[:]}$ and the set of all equivalence classes of families $((\sig_j)_{0\le j \le n}, (N_j)_{1\le j \le n})$,
where $n\geq 0$, $\sig_j$ are faces of $\sig$
such that $\sig_0=\{0\}$, $\sig_n=\sig$, $\sig_j\subsetneq  \sig_{j+1}$ for $0\leq j<n$, and $N_j$ is an element of the interior of $\sig_j/\sig_{j-1}:=(\sig_j+\sig_{j-1,\R})/\sig_{j-1,\R}$. Two such families 
$((\sig_j)_j, (N_j)_j)$ and $((\sig'_j)_j, (N'_j)_j)$ are equivalent if and only if $\sig'_j=\sig_j$ for $0\leq j\leq n$ and there are $c_j\in \R_{>0}$ such that $N_j'= c_jN_j$ for $1\leq j\leq n$.

In fact, for such a family $((\sig_j)_j, (N_j)_j)$, the corresponding $r\in \sig_{[:]}$ is given by  $r(f,g)= \tilde N_j(f)/\tilde N_j(g)$, where
 $j$ is the smallest integer such that at least one of $f$ and $g$ does not kill $\sig_j$ and $\tilde N_j$ is a lifting of $N_j$ to $\sig_j$.
 (\cite{KNU4} 4.1.6.)
  \end{para}
  
  \begin{para}\label{ofp} Let $(V, W, \sig, Y)$ be a monodromy system and let $p\in \sig_{[:]}$ which is the class of  $((\sig_j)_{0\leq j\leq n}, (N_j)_{1\leq j\leq n})$ as in \ref{1:1}. Let $\tilde N_j$ be a lifting of $N_j$ to $\sig_j$. Since $(\tilde N_1, \dots, \tilde N_n)$ is a Deligne system (Proposition \ref{p:ms&Ds}), we have the associated action of $\mathbb{G}_m \times \SL(2)^n$ (\ref{SL2}), $Y^j$ ($0\leq j\leq n$) and $\hat N_j$  ($1\leq j \leq n$).  $N_1$ and  $\hat N_j$ ($1\leq j\leq n$) are determined by $p$ modulo multiplications by elements of $\R_{>0}$ (to see that $\hat N_j$ is independent of the choice of the lifting $\tilde N_j$ of $N_j$, use the fact that the elements of $\sig_{j-1,\R}$ are of weight $-2$ for $Y^{j-1}$), and the action of $\mathbb{G}_m\times \SL(2)^n$ is determined by $p$ modulo the conjugate by $\begin{pmatrix} a_j^{-1}& 0\\0 & a_j\end{pmatrix}_{1\leq j\leq n}$ for $a_j\in \R_{>0}$. We have $Y^0=\spl_W(N_1+\sum_{j=2}^n \hat N_j)$ by Proposition \ref{twospl}.

  \end{para}
   
   \begin{para}\label{:Phi} If $\Phi$ is a set of faces of $\sig$ containing $\{0\}$ and $\sig$ and  totally ordered for the inclusion, and if  $\sig_{[:]}(\Phi)$ denotes the subset of $\sig_{[:]}$ consisting of elements whose corresponding pairs $((\sig_j)_{0\leq j\leq n}, (N_j)_{1\leq j\leq n})$ satisfy $\sig_j\in \Phi$ for all $0\leq j\leq n$, then $\sig_{[:]}(\Phi)$ is an open set of $\sig_{[:]}$.
   
   \end{para}
   
     \begin{para}\label{toric}
  For a sharp finitely generated cone $\sig$, let $\sig^{\circ}$ be the interior of $\sig$. Assume $\sig\neq \{0\}$. Then $\sig^{\circ}/\R_{>0}$ is regarded as a subset of $\sig_{[:]}$ because  it is regarded as the set of the equivalence classes in \ref{1:1} such that $n=1$. It is a dense open subset of $\sig_{[:]}$ and is a real analytic manifold. 
  
 In the rest of this section (Section \ref{s:ratio}),  we discuss the structure of $\sig_{[:]}$ as a real analytic manifold with corners in which corners appear outside $\sig^{\circ}/\R_{>0}$.
  \end{para}

\begin{para}\label{Njk} 

Let $\sig$ be a  sharp finitely generated cone. We endow $\sig_{[:]}$ a structure of a real analytic manifold with corners.

In the case $\sig=\{0\}$, $\sig_{[:]}$ is a one point set, and hence we have the structure evidently. Assume $\sig\neq\{0\}$.

By a {\it base along faces}, we mean a family $\Psi= ((\sig_j)_{0\leq j\leq n}, (N_{j,k})_{(j,k)\in J})$, where
$\sig_j$ are faces of $\sig$ such that $\{0\}=\sig_0 \subsetneq \sig_1 \subsetneq \cdots \subsetneq \sig_n=\sig$, $J:=\{(j,k)\;|\; 1\leq j\leq n, 1\leq k \leq r(j)\}$ with $r(j)=\dim_\R(\sig_{j,\R})-\dim_\R(\sig_{j-1, \R})$ 
and $(N_{j,k})_{(j,k)\in J}$ is an $\R$-base of $\sig_\R$ satisfying the following conditions (i) and (ii).

\smallskip

(i) For each $1\leq j\leq n$, $N_{j,k} \in \sig_j$ for $1\leq k\leq r(j)$ and $(N_{j,k}\bmod \sig_{j-1,\R})_{1\leq k\leq r(j)}$ is an $\R$-base of $\sig_{j,\R}/\sig_{j-1,\R}$.

(ii) For $(j,k)\in J$, the image of $N_{j,k}$ in $\sig_j/\sig_{j-1}$ (\ref{1:1}) belongs to the interior of $\sig_j/\sig_{j-1}$.

\smallskip

Let $U(\Psi)$ be the subset of $\sig_{[:]}$ consisting of all the classes of 
\begin{equation}\tag{\ref{Njk}.1}
\Bigl((\sig_{\varphi(i)})_{0\leq i\leq n'}, \bigl(\sum_{j=\varphi(i-1)+1}^{\varphi(i)}\sum_{k=1}^{r(j)} y_{j,k}N_{j,k}\bigr)_{1\leq i\leq n'}\Bigr)
\end{equation}
for some $n'$ such that $0\leq n'\leq n$,  some increasing injection $\varphi: \{0,\dots, n'\}\to \{0, \dots, n\}$ such that $\varphi(0)=0$ and $\varphi(n')=n$, and some
  $y_{j,k}>0$. 
  
 Then $U(\Psi)$ is an open set of $\sig_{[:]}$.
 
Further we can define an injective open map $$\iota_{\Psi}: U(\Psi)\to \R^{n-1}_{\geq 0} \times \R_{>0}^{J_1}, $$ 
where 
$J_1:=\{(j,k)\in J\;|\; k\neq 1\}$, which sends the class of (\ref{Njk}.1) in $U(\Psi)$  to $$((a_j)_{1\leq j\leq n-1}, (y_{j,k}/y_{j,1})_{(j,k)\in J_1}),$$
{where} $a_j=0$ if $j$ is in the image of $\varphi$, and $a_j= (y_{j+1, 1}/y_{j,1})^{1/2}$ otherwise.

Thus $$U(\Psi)\cap (\sig^{\circ}/\R_{>0})=\Bigl\{\sum_{(j,k)\in J} y_{j,k}N_{j,k}\;\Bigl|\; y_{j,k}>0\Bigr\}\bigl/\R_{>0},$$
and $\iota_{\Psi}$ sends the class of $\sum_{(j,k)\in J} y_{j,k}N_{j,k}\in \sig^{\circ}$ ($y_{j,k}>0$) to $$(((y_{j+1,1}/y_{j,1})^{1/2})_{1\leq j\leq n-1}, (y_{j,k}/y_{j,1})_{(j,k)\in J_1}).$$

For $p\in \sig_{[:]}$, we say that $p$ is {\it encased in $\Psi$} if $p$ is a class of 
 $((\sig_j)_{0\leq j\leq n}, (N_j)_{1\leq j\leq n})$ with $N_j=\sum_{k=1}^{r(j)} c_{j,k}N_{j, k}\bmod \sig_{j-1,\R}$ for some  $c_{j,k}\in \R_{>0}$, $c_{j,1}=1$,  for $1\leq j\leq n$. (We say that $p$ is encased in $\Psi$ with $(c_{j,k})_{j,k}$.) For such a $p$, 
$\iota_{\Psi}$ sends  $p$ to $((0)_{1\leq j \leq n-1}, (c_{j,k})_{(j,k)\in J_1})$. 

Each $p\in \sig_{[:]}$ is encased in some base $\Psi$ along faces. This 
is seen by applying to each $\sig_j/\sig_{j-1}$  the general fact that any open cone is covered by the open subcones generated by interior points.

Hence $U(\Psi)$ for bases $\Psi$ along faces form an open covering of $\sig_{[:]}$.

We endow  $\sig_{[:]}$ with the  unique structure of a real analytic manifold with corners of dimension $\dim_\R(\sig_\R)-1$ satisfying the following condition:

For every base $\Psi$ along faces, the map $\iota_{\Psi}$ is  an open immersion of real analytic manifolds with corners.

\end{para}

\begin{para}\label{Aan1} By \ref{Njk}, we have the sheaf of analytic functions on $\sig_{[:]}$. In the case where $x$ is an fs log point whose residue field is $\C$, this sheaf coincides with the sheaf of real analytic functions on $x_{[:]}$ denoted by $A_1^{\an}$ in \cite{KNU6} 1.4.19.

\end{para}

\begin{para}\label{narrow} There is a slightly different  structure of  a real analytic manifold with corners on $\sig_{[:]}$ which is characterized in the same way as in \ref{Njk}  except that we replace $(y_{j+1,1}/y_{j,1})^{1/2}$ in the definition of the immersion  $\iota_{\Psi}$  from $U(\Psi)$ by $y_{j+1, 1}/y_{j, 1}$. 

The sheaf of real analytic functions for this structure is called the sheaf of real analytic functions in the narrower sense. 

For example, if $\sig=\R_{\geq 0}^2$ as in Example \ref{eg}, we have an isomorphism $\sig_{[:]}\cong [0,\infty]\;;\; r \mapsto r(t_1, t_2)^{1/2}$  of real analytic manifolds with corners, and the homeomorphism $\sig_{[:]}\cong [0, \infty]$ in Example \ref{eg} without $(\,\cdot\,)^{1/2}$ is an isomorphism for the structure in the \lq\lq narrower'' sense. 

\end{para}

\begin{para}\label{mor[:]} 
For a homomorphism $h: \sig\to \sig'$ of sharp finitely generated cones, we have a canonical  map $\sig_{[:]}\to \sig'_{[:]}$ in the following two cases.

\medskip
(i) The image of $h$ contains some inner point of $\sig'$.

(ii) The cokernel of $\sig_\R\to \sig'_\R$ is of dimension $\leq 1$.  

\medskip
This map $\sig_{[:]}\to \sig'_{[:]}\;;\; r\mapsto r'$ is defined as follows. 
In the case (i), for homomorphisms $f,g:\sig'\to \R^{\add}_{\geq 0}$ such that $(f,g)\neq (0,0)$, we have $(f\circ h, g\circ h)\neq (0,0)$. We define $r'(f,g)= r(f\circ h, g\circ h)$. In the case where $\sig_\R\to \sig'_\R$ is surjective, the condition (i) is satisfied and we use the definition in the case (i). Assume that the cokernel of $\sig_\R\to \sig'_\R$ is of dimension one. We define $r'(f, g)= r(f,g)$ if $(f\circ h, g\circ h)\neq (0,0)$. If $(f\circ h, g\circ h)=(0,0)$  but $(f,g) \not=(0,0)$, by taking $N\in \sig'$ which is not in the image of $\sig_\R$, we define $r'(f,g)= f(N)/g(N)$, which is independent of the choice of $N$.

This map $\sig_{[:]}\to \sig'_{[:]}$ is  a morphism of real analytic manifolds with corners (also for the \lq\lq narrower sense'').

\end{para}

\begin{para}\label{[:]'} In \ref{mor[:]}, consider the case $\sig'=\sig \times \R_{\geq 0}$ and  $h(a)=(a,0)$ for $a\in \sig$. For  $p\in \sig_{[:]}$ 
corresponding to $((\sig_j)_{0\leq j\leq n}, (N_j)_{1\leq j\leq n})$, the image of $p$ in $\sig'_{[:]}$ is the class of $((\sig'_j)_{0\leq j\leq n+1}, (N'_j)_{1\leq j\leq n+1})$, where $\sig'_j=\sig_j$ and $N'_j=N_j$ for $1\leq j\leq n$, $\sig'_{n+1}=\sig'$, and $N'_{n+1}=(0,1)$. 
The map $\sig_{[:]} \to (\sig')_{[:]}$ is injective and we  will identify $\sig_{[:]}$ with its image in $(\sig')_{[:]}$. 

Note that we have an isomorphism $\sig^{\circ}\overset{\cong}\to (\sig')^{\circ}/\R_{>0}$ by $N\mapsto (N, 1)$. 
 
There is a strong relation (Proposition \ref{sigtor}) between $(\sig')_{[:]}$ and a toric variety. Assume $\sig=\Hom(\cS, \R^{\add}_{\geq 0})$ for an fs monoid $\cS$.
We consider the diagram with a cartesian square
\begin{equation*}\tag{$*$} 
\begin{matrix}  |\torus|(\cS)_{<1}  & \subset & |\toric|(\cS)_{[:],<1} & \supset & 0_{[:]}\\
\cap && \cap && \\
|\torus|(\cS) &\subset & |\toric|(\cS)_{[:]} & \subset & \toric(\cS)_{[:]}.\end{matrix}
\end{equation*}
Here $\toric(\cS)_{[:]}$ is the space of ratios associated to  the toric variety $\toric(\cS)=\Hom(\cS, \C^{\mult})$ ($\C^{\mult}=\C$ regarded as a multiplicative monoid), $|\toric|(\cS)_{[:]}$ is the inverse image of $|\toric|(\cS)= \Hom(\cS, \R_{\geq 0}^{\mult})$ under $\toric(\cS)_{[:]} \to \toric(\cS)$, $|\torus|(\cS)=\Hom(\cS, \R_{>0})$. Let $|\toric|(\cS)_{<1}$ (resp.\ $|\torus|(\cS)_{<1}$) be the set of homomorphisms $\cS\to \R^{\mult}_{\geq 0}$ (resp.\ $\cS\to \R_{>0}$) which send $\cS\smallsetminus \{1\}$ to $\{x\in \R\;|\;x<1\}$. Let $|\toric|(\cS)_{[:],<1}$ be the inverse image of $|\toric|(\cS)_{<1}$ under $|\toric|(\cS)_{[:]} \to |\toric|(\cS)$. 
Let $0\in \toric(\cS)$ be the fs log point corresponding to the homomorphism $\cS\to \C^{\mult}$ which sends all $\cS\smallsetminus \{1\}$ to $0$.

  Then $|\toric|(\cS)_{[:],<1}$ is  an open neighborhood of $0_{[:]}$ in $|\toric|(\cS)_{[:]}$.
  
   For example, in the case $\cS=\N^2$, the above diagram becomes
  $$\begin{matrix}  (0,1)^2  & \subset & ([0, 1)^2)_{[:]}  & \supset & 0_{[:]}\\
  \cap && \cap &&\\
  \R^2_{>0}& \subset& (\R^2_{\geq 0})_{[:]} &\subset & (\C^2)_{[:]} \end{matrix}$$
  in which every space of the form  $(\,\cdot\,)_{[:]}$ is obtained from the space $(\,\cdot\,)$ by replacing the point $(0,0)$ in $(\,\cdot\,)$ by $0_{[:]}=[0, \infty]$. 
A directed family $(a_{\la}, b_{\la})$ in $\R_{\geq 0}^2\smallsetminus \{(0,0)\}$ converges in $(\R^2_{\geq 0})_{[:]}$ to 
 $r\in 0_{[:]}=[0,\infty]$ if and only if $(a_{\la}, b_{\la})$ converges to $(0,0)$ in $\R^2_{\geq 0}$ and  $a_{\la}/b_{\la}\in [0,\infty]$ converges to $r$ in $[0,\infty]$.

On the other hand, let $(\sig')_{[:],<1}$ be the set of points of $(\sig')_{[:]}$ which are classes of $((\sig'_j)_{0\leq j\leq n}, (N'_j)_{1\leq j\leq n})$ such that $\sig'_j$ for $0\leq j<n$ are faces of $\sig$. By \ref{:Phi}, this is an open neighborhood of $\sig_{[:]}$ in $(\sig')_{[:]}$.

\end{para} 

\begin{prop}\label{sigtor} The part of the above diagram $(*)$ consisting of the upper row and the left column is canonically isomorphic to
$$\begin{matrix}  \sig^{\circ} &\subset& (\sig')_{[:],<1} & \supset & \sig_{[:]}\\
\cap && && \\
\sig_\R.&&&&
\end{matrix}$$
as a diagram of topological spaces.
\end{prop}

\begin{pf} We have the isomorphism $|\torus|(\cS)=\Hom(\cS, \R_{>0})\cong \sig_\R=\Hom(\cS, \R)$ induced by $\R_{>0}\overset{\cong}\to \R\;;\;x\mapsto -(2\pi)^{-1}\log(x)$. 
This induces $|\torus|(\cS)_{<1}\cong \sig^{\circ}$. 

We consider $|\toric(\cS)|_{[:], <1}$. An element $p: \cS\to \R^{\mult}_{\geq 0}$ of $|\toric|(\cS)_{<1}$ is understood as a pair $(\tau, N)$, where $\tau$ is a face of $\sig$ and $N$ is an element of the interior of the cone $\sig/\tau= (\sig +\tau_\R)/\tau_\R$, as follows. $\tau$ is the face of $\sig$ corresponding to the face $\cT= p^{-1}(\R_{>0})\subset \cS$ of $\cS$. $N$ is the element corresponding to the homomorphism $\cT \to \R^{\add}_{\geq 0}\;;\; x\mapsto -(2\pi)^{-1}\log(p(x))$. The monodromy cone $\sig(p)$ at $p\in \toric(\cS)$ is $\tau$. Hence $|\toric|(\cS)_{[:],<1}$ is identified with the triple $(\tau, N, a)$, where $(\tau,N)$ is as above and $a$ is an element of $\tau_{[:]}$. Such a triple bijectively corresponds to  an element of $(\sig')_{[:],<1}$: If $a$ is the class of $((\tau_j)_{0\leq j\leq n}, (N_j)_{1\leq j\leq n})$, we have the class of  $((\sig'_j)_{0\leq j\leq n+1}, (N'_j)_{1\leq j\leq n+1})$ in $(\sig')_{[:],<1}$, where $\sig'_j=\tau_j$ and $N'_j=N_j$ for $1\leq j\leq n$, $\sig'_{n+1}=\sig'$, and $N_{n+1}=(N, 1)$.
\end{pf} 

\begin{para}\label{U'}

Let the notation be as in  \ref{[:]'}.  Assume $\sig\neq\{0\}$. 
Then from  a base 
$\Psi=((\sig_j)_{0\leq j\leq n}, (N_{j,k})_{(j,k)\in J})$ along faces 
for $\sig$ in \ref{Njk}, we obtain a base  $\Psi'=(\sig'_{0 \leq j\leq n+1}, (N'_{j,k})_{(j,k)\in J'})$  along faces for $\sig'$ as follows: $\sig'_j=\sig_j$ for $0\leq j\leq n$, $\sig'_{n+1}=\sig'$, $J'= J\;\text{\scriptsize $\coprod$}\;\{(n+1, 1)\}$, $N'_{j,k}=N_{j,k}$ if $j\neq n+1$, and $N'_{n+1,1}= (0,1)\in\sig\times \R_{\geq 0}$. 

Then for the map $(\sig')_{[:]}\to \sig_{[:]}$ induced by the projection $\sig'\to \sig\;;\;(a, 0)\mapsto a$ (\ref{mor[:]}), $U(\Psi')$ coincides with the inverse image of $U(\Psi)$ in $(\sig')_{[:],<1}$. Hence these $U(\Psi')$ form an open covering of $(\sig')_{[:],<1}$. 

We have the commutative diagram
$$\begin{matrix}  U(\Psi') &\overset{\iota_{\Psi'}}\to & \R^n_{\geq 0} \times \R_{>0}^{J_1} \\
\downarrow && \downarrow \\
U(\Psi) & \overset{\iota_{\Psi}}\to & \R^{n-1}_{\geq 0} \times \R_{>0}^{J_1},\end{matrix}$$
where $\R^n_{\geq 0}\to \R^{n-1}_{\geq 0}$ is the projection to forget the $n$-th factor. We have $U(\Psi')\cap \sig_{[:]}=U(\Psi)$, and $U(\Psi)$ is the subset of $U(\Psi')$ consisting of all elements whose $n$-th factors in $\R^n_{\geq 0}$  of the images by $\iota_{\Psi'}$ are $0$. 

From this, we see that $(\sig')_{[:],<1} \smallsetminus \sig_{[:]}\to \sig_{[:]}$ is an $\R_{>0}$-torsor in the category of real manifolds with corners (also for the \lq\lq narrower''  sense).

$(U(\Psi')\cap (\sig')^{\circ})/\R_{>0}$ is identified, via $(\sig')^{\circ}/\R_{>0} \cong (\sig)^{\circ}$,  with the subset of $\sig^{\circ}$ consisting of elements $\sum_{(j,k)\in J} y_{j,k}N_{j,k}$ with $y_{j,k}>0$.

\end{para}

\begin{para}\label{Hratio}

In the theory of degeneration of Hodge structure, the notion nilpotent orbit is important. In the formulation in \cite{KNU3} Section 2.2, a nilpotent orbit (we always assume that its  $\gr^W$ is polarizable) is 
$(H_\R, W, \sig, F)$, where $H_\R$ is a finite dimensional $\R$-vector space, $W$ is a finite increasing filtration on $H_\R$, $\sig$ is a finitely generated cone with an admissible action on $H_\R$, and $F$ is a finite decreasing filtration on $H_\C=\C\otimes_\R H_\R$ (called the Hodge filtration) satisfying certain conditions. 
One condition is that $(W(\sig),F)$ is a mixed Hodge structure. If $Y$ denotes the canonical splitting of $W(\sig)$  with respect to $F$ (in the sense of \cite{KNU08} Section 1),  we have a monodromy system $(H_\C, W_\C, \sig, Y)$. 

   Thus we can say that a monodromy system is like a nilpotent orbit without Hodge filtration and hence should be important in Hodge theory. But in this paper, we are trying to apply monodromy systems to non-archimedean geometry (as in Sections \ref{s:na} and \ref{s:reg})  and also to the logarithmic geometry over a finite field.

Assume  $\sig=\Hom(\cS, \R^{\add}_{\geq 0})$ for a sharp fs monoid $\cS$. So,  the above  nilpotent orbit is regarded as a log mixed Hodge structure on the fs log point $0\in \toric(\cS)$ in \ref{[:]'}. 

In the theory of nilpotent orbit, the behavior of the map 
$$\sig_\R \ni N \mapsto \exp(iN)F$$
 \lq\lq at the limit'' is important. In the theory of SL(2)-orbit, the meaning of this \lq\lq at the limit'' is understood in \cite{KNU6} 4.5.5, 4.5.6 as what happens  when a point of $\sig_\R$ approaches  a point of $0_{[:]}$ in $|\toric|(\cS)_{[:]}$ (\ref{[:]'}). 
We reinterpret this as what happens when a point of $\sig^{\circ}$  approaches a point of $\sig_{[:]}$ in $(\sig')_{[:]}$ (Proposition \ref{sigtor}), and we use this idea in Sections \ref{s:spl} and \ref{s:del} when we consider Hodge theory in \ref{Hodge1} and \ref{delH3}. 

\end{para}

\section{Analogy}\label{s:na}
  We explain more the analogy mentioned in Introduction.

\begin{para}\label{analogy} 
We have the strong analogy between the following (i) and (ii).

\medskip

(i) $S=\Delta^n$ with $\Delta$ being the unit disc over $\C$. We endow $S$ with the log structure generated by coordinate functions $t_1, \dots, t_n$. Let  $U=(\Delta^*)^n$, where $\Delta^*=\Delta\smallsetminus \{0\}$. 

(ii) $S=\Spec(O_K[[t_1, \dots, t_n]])$ for a complete discrete valuation field $K$ with finite residue field $k$ and with valuation ring $O_K$. We  endow $S$  with the log structure generated by $t_1, \dots, t_n, t_{n+1}$, where $t_{n+1}$ denotes  a prime element of $O_K$. Let $U= \Spec(O_K[[t_1, \dots,t_n]][t_1^{-1}, \dots, t_{n+1}^{-1}])$.

\medskip

In (i), let $0:=(0, \dots, 0)\in S$. 
In (ii), let $0$ be the unique closed point of $S$. 
Both in (i) and (ii), we endow $0$ with the log structure induced from that of $S$. 

In (i), consider a variation of mixed Hodge structure $H$ on $U$. It often extends to a log mixed Hodge structure on $S$ and induces at $0$ a nilpotent orbit $(N_1, \dots, N_n, F)$. The asymptotic behavior of $H$ on  $U$ at points near $0$ is understood well by this nilpotent orbit.

In (ii), let $\ell$ be a prime number which is different from the characteristic of $k$, and consider a smooth $\bar \Q_{\ell}$-sheaf $H$ on the \'etale site of  $U$. It often extends to a smooth $\bar \Q_{\ell}$-sheaf on the log \'etale site of $S$ and induces at $0$ an object of the category $\cA_0$ which gives a monodromy system $(N_1, \dots, N_{n+1})$. The asymptotic behavior of $H$ on  $U$ at points near $0$ is understood well by this monodromy system. 

Here in (i),  $(N_j)_{1\leq j\leq n}$ is the dual base of the base $(t_j)_{1\leq j\leq n}$ of $(M_S/\cO_S^\times)_0$. 
In (ii), $(N_j)_{1\leq j\leq n+1}$ is  the dual base of $(t_j)_{1\leq j\leq n+1}$ of $(M_S/\cO_S^\times)_{\bar 0}$.

In the theory of nilpotent orbit $(N_1, \dots, N_n, F)$, it is important to consider $\exp(iy_1N_1+\dots+iy_nN_n)F$ ($y_j>0$). In the theory of monodromy system $(N_1, \dots, N_{n+1})$, it is important to consider  $y_1N_1+\dots+y_nN_n+N_{n+1}$ ($y_j>0$). See \ref{Natp} for a basic fact about the latter.
\end{para}

\begin{para}\label{Ndef} We review the monodromy operator of a Galois representation of a local field. 

 Let $K$ be a complete discrete valuation field  with residue field $k$. Let $\ell$ be a prime number which is different from the characteristic of $k$.

 Let $V$ be a finite dimensional $\Q_{\ell}$-vector space with a continuous action of $G_K=\Gal(\bar K/K)$, where $\bar K$ is the separable closure of $K$. Assume 
that the action of the inertia subgroup $I$ of $G_K$ on $V$ is quasi-unipotent. (This condition is satisfied if $V=H^m(X)_{\ell}:=H^m_{\et}(X \otimes_K \bar K, \Q_{\ell})$ for a scheme $X$ of finite type over $K$.) Then we have a nilpotent linear map  $$N: V\to V(-1)$$
(called the {\it monodromy operator}) defined  as follows.
Take a prime element $\pi$ of $K$. Then we have the non-trivial homomorphism $$h: I\to \Z_{\ell}(1)\;;\; \sig\mapsto (\sig(\pi^{1/\ell^n})/\pi^{1/\ell^n})_{n \geq 1}$$ which is independent of the choice of $\pi$.  
   Let 
$$N=h(\sig)^{-1}\log(\sig): V\to V(-1),$$ 
where $\sig$ is an element of  $I$  such that $h(\sig)\neq 0$ and such that the action of $\sig$ on $V$ is unipotent.
Then $N$ is independent of the choice of such a $\sig$.

\end{para}

\begin{para}\label{Ndef2} Let $K$ and $V$ be as in \ref{Ndef}. Assume that $K$ is a finite extension of a complete discrete valuation field $J$ (we assume that the valuation of $K$ comes from  that of $J$). Then the 
{\it monodromy operator of $K$ on $V$ with respect to $J$} is defined to be $e(K/J)^{-1}N$, where $N$ is the monodromy operator in \ref{Ndef} and $e(K/J)$ is the ramification index.

If $K'$ is a finite extension of $K$, the monodromy operator of $K'$ on $V$ with respect to $J$ coincides with the monodromy operator of $K$ on $V$ with respect to $J$.

 \end{para}

\begin{para}\label{Natp} Consider the  case (ii) in \ref{analogy}. 
Let $p$ be a closed point of $U$, and let $K(p)$ be the residue field of $p$. Then the monodromy operator of the local field $K(p)$ on the stalk of $H$ with respect to $K$ (\ref{Ndef2}) is equal to 
$$\bigl(\sum_{j=1}^n y_jN_j\bigr)+N_{n+1} \quad \text{with}  \; y_j:=v(t_j(p)),$$ 
where $v$ is the additive valuation $K(p)^\times \to \Q$ which sends a prime element of $K$ to $1$. 

\end{para}

In \ref{ht1}--\ref{ht5}, as an example of \ref{analogy}, we compare the archimedean theory and the non-archimedean theory of the local height pairing on an elliptic curve.  
In \ref{ht3}--\ref{ht5},  we take $n=1$ in \ref{analogy}, and take a family of elliptic curves $E$ over $U$ which degenerates at $0$, and consider $H$ in \ref{analogy}  obtained from a pair $(A,B)$ of  families of $0$-cycles on $E$ with degree $0$. There we consider the  asymptotic behavior of the local height pairing of $A$ and $B$ as the asymptotic behavior of the invariant $\delta_W$ of $H$.

In \ref{ht1}, we review the local height pairing. In \ref{ht2}, we understand it as $\delta_W$. In \ref{ht3}--\ref{ht5}, we consider the asymptotic behavior. 

\begin{para}\label{ht1} Classical theory.

Let $K$ be either 

(i) $\C$ or 

(ii) a complete discrete valuation field with finite residue field.

In the case (ii), let $v:K^\times \to \Z$ be the normalized additive valuation of $K$. 

Let $E$ be an elliptic curve over $K$. 
In the case (ii), assume that $E$ is a Tate curve. 
We have $E(K)= K^\times/q^{\Z}$ for some $q\in K^\times$ such that $|q|<1$.

 Let $$\theta_q(u)= \prod_{n=0}^{\infty} (1-q^nu) \cdot \prod_{n=1}^{\infty} (1-q^nu^{-1})\in \Z[u^{\pm 1}][[q]].$$ We have $\theta_q(qu)= -u^{-1}\theta_q(u)$.

Consider two divisors $A$ and $B$ on $E$ of degree $0$:

\noindent $A= \sum_j m(j)(\alpha_j\bmod q^{\Z})$, $B=\sum_j n(j)(\beta_j\bmod q^{\Z})$, $\alpha_j,\beta_j\in K^\times$, $\sum_j m(j) =\sum_j n(j)=0$. Assume that the supports of $A$ and $B$ are disjoint.

The  local height pairing $\langle A, B\rangle$ is given by 
$$\langle A, B\rangle= (2\pi)^{-1}\Bigl(\sum_{j,h} m(j)n(h)\log|\theta_q(\alpha_j/\beta_h)|- \frac{(\sum_j m(j)\log|\alpha_j|) \cdot (\sum_j n(j)\log|\beta_j|)}{\log|q|}\Bigr)$$
in the case (i), and 
$$\langle A, B\rangle= \sum_{j,h} m(j)n(h)v(\theta_q(\alpha_j/\beta_h))- \frac{(\sum_j m(j)v(\alpha_j)) \cdot (\sum_j n(j)v(\beta_j))}{v(q)}$$
in the case (ii).

\end{para}

\begin{para}\label{ht2} Mixed Hodge structure in the case (i) and monodromy system with one $N$ in the case (ii).

We have an interpretation of $\langle A, B\rangle$ as 
$$\langle A, B\rangle =\delta_W(H_{A,B}).$$ Here $H_{A,B}$ is a mixed Hodge structure  in the case (i) and a monodromy system  with one $N$ in the case (ii) defined as follows.

It is of rank $4$ with base $e_j$ ($j=1,2,3,4$). $W_{-3}=0$, $W_{-2}$ is generated by $e_1$, $W_{-1}$ is generated by  $e_1, e_2, e_3$, and $W_0$ is the total space. 

In the case (i) (see \cite{KNU2} 4.4): $F^{-1}=H$, $F^1=0$, $F^0$ is generated by $e_3+\tau e_2+we_1$ and $e_4+ze_2+\la e_1$, where $\tau=(2\pi i)^{-1}\log(q)$, $z=(2\pi i)^{-1}\sum_j m(j)\log(\alpha_j)$, $w=(2\pi i)^{-1}\sum_j n(j)\log(\beta_j)$, $\la=(2\pi i)^{-1}\sum_{j,h} m(j)n(h)\log\theta_q(\alpha_j/\beta_h)$.

In the case (ii): $N(e_1)=N(e_2)=0$. $N(e_3)= v(q) e_2 + (\sum_j n(j) v(\beta_j))e_1$, $N(e_4)= (\sum_j m(j)v(\alpha_j))e_2+ (\sum_{j,h} m(j)n(h)v(\theta_q(\alpha_j/\beta_h)))e_1$.  $Y(e_4)=0$, $Y(e_3)=0$, $Y(e_2)= -2e_2$, $Y(e_1)=- 2e_1$. This is an object of $\cA_0$ for $0=\Spec(k)$, where  $0$ is endowed with the inverse image of the canonical log structure of $\Spec(O_K)$, obtained from a representation of $\Gal(\bar K/K)$ with a unipotent action of the inertia group.

\end{para}

\begin{para}\label{ht3} Variation and degeneration. 

In the case (i), let $\Delta$ be the  unit disk with the coordinate function $t$ and let $U=\Delta\smallsetminus \{0\}$. 

In the case (ii), let $\Delta:=\Spec(O_K[[t]])$ and let $U=\Spec(O_K[[t]][1/t,1/\pi])$, where $\pi$ is a prime element of $K$.

 In each case (i), (ii), we assume the following. 
 
  Case (i).  Let $q\in t^c\cO(\Delta)^\times$ for some integer $c>0$ and  $\alpha_j\in t^{a(j)}\cO(\Delta)^\times$, $\beta_j\in t^{b(j)}\cO(\Delta)^\times$ for some $a(j), b(j)\in \Z$. We assume $\alpha_j -q^r\beta_h$ has no zero on $U$ for every $j,h$ and every $r\in \Z$.

  This  implies 
  $\theta_q(\alpha_j/\beta_h)\in t^{d(j,h)}\cO(\Delta)^\times$ for some $d(j,h)\in \Z$. 

Case (ii). Let $q\in t^c\pi^{c'}\cO(\Delta)^\times$ for some integers $c>0,c'>0$, and $\alpha_j\in t^{a(j)}\pi^{a'(j)}\cO(\Delta)^\times$, $\beta_j\in t^{b(j)}\pi^{b'(j)}\cO(\Delta)^\times$ for some $a(j), a'(j), b(j), b'(j)\in \Z$. We assume $\alpha_j-q^r \beta_h \in t^{\Z}\pi^{\Z}\cO(\Delta)^\times$ for every $j,h$ and every $r\in \Z$. This  implies 
 $\theta_q(\alpha_j/\beta_h) \in t^{d(j,h)}\pi^{d'(j,h)}\cO(\Delta)^\times$ for some $d(j,h), d'(j,h)\in \Z$. 

Let $m(j), n(j)\in \Z$, $\sum_j m(j)=\sum_j n(j)=0$. 

We have the variation $(E_{q(p)})_p$ of elliptic curves, where $p$ ranges over $K$-points of $U$,  and the family of divisors $A= \sum_j m(j)(\alpha_j)$, $B=\sum_j n(j)(\beta_j)$ of degree $0$ on it. 

In the case (i), we have the variation $H_{A,B}$ of mixed Hodge structure on $U$ whose restriction to each $p\in U$ gives the mixed Hodge structure in \ref{ht2}. In the case (ii), for a prime number $\ell$ which is not the characteristic of $k$, we have a smooth $\Q_{\ell}$-sheaf $H_{A,B}$ on $U$ whose restriction to each $p\in U(K)$ gives the monodromy system in \ref{ht2}. 

\end{para}

\begin{para}\label{ht4} Nilpotent orbit in the case (i) and a monodromy system with two $N$ in the case (ii).

In the case (i) (resp.\ (ii)), 
endow $\Delta$ with the log structure generated by $t$ (resp.\ $t$ and $\pi$). Then we have the variation of log mixed Hodge structure $H_{A,B}$ on $\Delta$  (resp.\ a $\Q_{\ell}$-sheaf $H_{A,B}$ on the log \'etale site of  $\Delta$). Its restriction $H_{A,B}(0)$ to the fs log point $\Spec(\C)$ (resp.\ $\Spec(k)$) at $t=0$ endowed with the log structure induced by that of $\Delta$ is identified with the following nilpotent orbit $(N_1, F)$ in (i) (resp.\ monodromy system $(N_1, N_2, Y)$ in (ii)).

 In both (i) and (ii),  $N_1(e_4)= ae_2+de_1$, $N_1(e_3)= ce_2+be_1$, $N_1(e_2)=N_1(e_1)=0$,  where $$a=\sum_j m(j)a(j), \quad   b=\sum_j n(j)b(j)
\quad d=\sum_{j,h}  m(j)n(h)d(j,h).$$
 
In  (i),    $F$ is like $F$ in \ref{ht2} but we replace $\log(X)$ for $X=q, \alpha_j, \beta_j, \theta_q(\alpha_j/\beta_h)$ there by $\log(Y)$, where $Y$ is the value of $t^{-c}q$, $t^{-a(j)}\alpha_j$, $t^{-b(j)}\beta_j$, $t^{-d(j,h)}\theta_q(\alpha_j/\beta_h)$, at $0\in \Delta$, respectively.
 
In (ii), $N_2(e_4)=a'e_2+ d'e_1$, $N_2e_3= c'e_2+b'e_1$, $N_2(e_2)=N_2(e_1)=0$, where $a'= \sum_j m(j)a'(j)$,  $b'= \sum_j n(j)b'(j)$, $d'=\sum_{j,h} m(j)n(h)d'(j,h)$. $Ye_4=0$, $Ye_3= 0$, $Ye_2=-2e_2$, $Ye_1=-2e_1$.
  \end{para}

\begin{para}\label{ht5} Asymptotic behavior. 

In the case (i) (resp.\ (ii)), consider $p\in U$ (resp.\ $p\in U(K)$) with  $t(p)\to 0$, that is, $y\to \infty$, where $y=-(2\pi)^{-1}\log|t(p)|$ (resp.\ $y=v(t(p))$). 

Then both $\langle A(p), B(p)\rangle=\delta_W(H_{A,B}(p))$ and $\delta_W(\exp(iyN_1)F)$  have (resp.\ $\langle A(p), B(p)\rangle= \delta(yN_1+N_2)$ (\ref{Natp}) has) the asymptotic behavior $$\Bigl(d- \frac{ab}{c}\Bigr)y+O(1),$$
 where $a,b,d$ are as in \ref{ht3}.

\end{para}

\begin{para}\label{Pearl}
A general result on the asymptotic behavior of  the local height pairing of algebraic cycles in the archimedean case was obtained in Pearlstein \cite{Pe}. Its  non-archimedean version is Theorem \ref{thmht}. (The weight-monodromy  conjecture in the assumption of Theorem \ref{thmht} is true in the present  situation (ii)).

\end{para}
 
\section{SL(2)-orbit theorems on the torus twist back}\label{s:twi}

\begin{para}\label{tw1} Let $(V, W,\sig, Y)$ be a monodromy system with $\sig \neq \{0\}$.

Let $\Psi=((\sig_j)_{0\leq j\leq n}, (N_{j,k})_{(j,k)\in J})$ be 
a base along faces (\ref{Njk}).

Take a homomorphism $(\tau_j)_{1\leq j\leq n-1}: \mathbb{G}_m^{n-1} \to \Aut(V)$ of algebraic groups over $\C$ such that $\tau_j$ is a splitting of $W(\sig_j)$ for $1\leq j\leq n-1$ and such that for $1\leq j\leq k\leq n-1$, all elements of $\sig_j$ are of weight $-2$ for $\tau_k$.

For $y=(y_{j,k})_{(j,k)\in J}\in \R_{>0}^J$, 
let
$$t(y):= \prod_{j=1}^{n-1} \tau_j((y_{j+1,1}/y_{j,1})^{1/2}).$$ 

Let $R$ be the ring of polynomials over $\R$ in $(y_{j+1,1}/y_{j,1})^{1/2}$ ($1\leq j\leq n-1$) and $y_{j,k}/y_{j,1}$ ($(j,k)\in J_1$). 

\end{para}

The following theorem is a refined version of \cite{KNU} Theorem 5.3 and the proof given below is by the method there.

\begin{thm}\label{cRNy}  

Let the situation be as in $\ref{tw1}$. 
Let 
$$N_y:=t(y)^{-1}\sum_{(j,k)\in J} y_{j,k}N_{j,k}t(y)$$ with $y_{j,k}\in \R_{>0}$, $y_{n,1}=1$.   
  Then the following holds. 

$(1)$ $N_y \in R\otimes_{\R} \End(V)$. 

$(2)$ Let $p\in \sig_{[:]}$ and assume that $p$ is the class of $((\sig_j)_{0\leq j\leq n}, (N_j)_{1\leq j\leq n})$, $N_j=\sum_{k=1}^{r(j)} c_{j,k}N_{j,k}\bmod \sig_{j-1,\R}$ ($c_{j,k}\in \R_{>0}$, $c_{j,1}=1$) for $1\leq j\leq n$. In $\ref{tw1}$,take $(\tau_j)_j$ which corresponds to the splittings $Y^j$ of $W(\sig_j)$ associated to $p$ ($\ref{ofp}$). Then the image of $N_y$ under $$R\to \R\;\;;\;\; (y_{j+1,1}/y_{j,1})^{1/2}\mapsto 0, \;\; y_{j,k}/y_{j,1}\mapsto c_{j,k}$$ is equal to  $N_1+\sum_{j=2}^n \hat N_j\in \End(V)$.

\end{thm}

\begin{pf}
(1)
Since  $N_{j,k}$ is of weight $-2$ for $\tau_{j'}$ ($j'\geq j$), we have 
\begin{equation*}\tag{\ref{cRNy}.1}
t(y)^{-1}y_{j,k}N_{j,k}t(y)= y_{j,k}/y_{j,1} \cdot t_{<j}(y)^{-1}N_{j,k}t_{<j}(y),
\end{equation*}
where
$t_{<j}(y)= \prod_{j'=1}^{j-1} \tau_{j'}((y_{j'+1,1}/y_{j',1})^{1/2})$. For $1\leq j'< j$, since $N_{j,k}$ respects $W^{j'}$,  the $\tau_{j'}$-weights of $N_{j,k}$ are $\leq 0$. Hence $t_{<j}(y)^{-1}N_{j,k}t_{<j}(y)$ is a polynomial in $(y_{j'+1,1}/y_{j',1})^{1/2}$ for $1\leq j'\leq j-1$.

(2) follows from (\ref{cRNy}.1).
\end{pf}

\begin{para}\label{Njk2} In several theorems in this paper, we will consider the following situation:

$\sig$ is a sharp finitely generated cone. $\sig\neq \{0\}$.

$\Psi=((\sig_j)_{0\leq j\leq n}, (N_{j,k})_{(j,k)\in J})$ is a base along faces, $p\in \sig_{[:]}$ and $p$ is encased in $\Psi$ with $c_{j,k}$ (\ref{Njk}).  

$y_{j,k}>0$ for $(j,k)\in J$. 

We assume that the element
$\sum_{(j,k)\in J} y_{j,k}N_{j,k} \mod \R_{>0}$ of $\sig^{\circ}/\R_{>0}$ is sufficiently near to $p$. In other words, $y_{j+1, 1}/y_{j,1}$ ($1\leq j\leq n-1$) and $y_{j,k}/y_{j,1}-c_{j,k}$ ($(j,k)\in J_1$) are sufficiently near to $0$.

\end{para}

\begin{thm}\label{delthm2} Let $(V, W, \sig, Y)$ be a monodromy system. Let the situation be as in $\ref{Njk2}$. Let  $Y^j$ be the splitting of $W(\sig_j)$ associated to $p$ ($\ref{ofp}$). Let $N_y$ be as in Theorem $\ref{cRNy}$ which we define using $t(y)$ defined by $\tau_j$ associated to $p$. Then we have the following.

$(1)$  $\spl_W(N_y) = u^{\natural}(y)Y^0u^{\natural}(y)^{-1}$, 
where $u^{\natural}(y)$ is a convergent Taylor series in $(y_{j+1,1}/y_{j,1})^{1/2}$ ($1\leq j\leq n-1$) and  $y_{j,k}/y_{j,1}  -c_{j,k}$ ($(j,k)\in J_1$)  with coefficients in $\End(V)$ such that $u(y)-1\in W_{-1}\End(V)$ and such that $u(y)$ is of $Y^n$-weight $0$.

$(2)$  $\delta_W(N_y)$ with $y_{n,1}=1$ is a convergent Taylor series in 
 $(y_{j+1,1}/y_{j,1})^{1/2}$ ($1\leq j \leq n-1$) and $y_{j,k}/y_{j,1}  -c_{j,k}$ ($(j,k)\in J_1$) with coefficients in $\gr^W\End(V)$.

\end{thm}

\begin{pf} This follows from Theorem \ref{cRNy} by the fact $N_1+\sum_{j=2}^n \hat N_j\in \cD$ for $W'= W(\sig)$ (\ref{calD}, Proposition \ref{MSD2}) and by the analyticity of the maps $\spl_W$ and $\delta_W$ on $\cD$ (\ref{calD}). 
 \end{pf}

\begin{para} 
An analogue of Theorem \ref{cRNy} in Hodge theory  with pure $W$ is 
 \cite{CKS} Theorem 4.20 (viii) which says 
that $t(y)^{-1}\exp(\sum_{j=1}^n iy_jN_j)F$ (the Hodge analogue of $N_y$) with $t(y) =\sum_{j=1}^n \exp(-\frac{1}{2}\sum_{j=1}^n \log(y_j)Y_j)$ is expressed as ${}^eg(y)\exp(\sum_{j=1}^n i\hat N_j)\hat F_n$. A Hodge analogue with mixed $W$ is \cite{KNU2} Theorem 2.4.2 (ii).

A Hodge analogue of Theorem \ref{delthm2} (2) is \cite{KNU08} Theorem 0.5 (4).

\end{para}

\section{SL(2)-orbit theorems in one variable}\label{s:1var}

In this section, for a monodromy system $(N_1, N_2)$, we prove Theorems \ref{thm1var} and \ref{mild1var} on the asymptotic behaviors of $\spl_W(yN_1+N_2)$ and $\delta_W(yN_1+N_2)$ for $y\gg 0$. The results for  general monodromy systems are proved in Sections \ref{s:spl} and \ref{s:del} by the reduction to the results in this section.

\begin{thm}\label{thm1var} Let $(N_1, N_2)$ be a monodromy system. Let $W^j$ and $Y^j$ be those associated to  the Deligne system $(N_1,N_2)$.

$(1)$   For $y\gg 0$, we have $$\spl_W(yN_1+N_2)= u(y)Y^0u(y)^{-1},$$ where $u(y)=1+\sum_{m=1}^{\infty} u_my^{-m}$ is a convergent Taylor series in $y^{-1}$ with coefficients $u_m\in W_{-1}\End(V)$. For every $m\geq 1$, we have that $u_m\in W^1_{m-1}\End(V)$  and  $u_m$ is of $Y^2$-weight $0$.

$(2)$ For $y\gg 0$, we have $$\delta_W(yN_1+N_2)= \sum_{m=-1}^{\infty} \delta_m y^{-m},$$  where $\delta_m\in W_{-2}\gr^W\End(V)$ and $\sum_{m=0}^{\infty} \delta_my^{-m}$ is a convergent Taylor series in $y^{-1}$. We have $\delta_{-1}= \delta_W(N_1)$. For every $m\geq -1$, we have that  $\delta_m\in  W^1_{m-1}\gr^W\End(V)$ and  $\delta_m$ is of $Y^2$-weight $-2$.

\end{thm}

\begin{para} Let $\C(y)$ be the rational function field in one variable $\C(y)$. By the theory of Deligne systems considered over the field $\C(y)$, we have that $\spl_W(yN_1+N_2)=u(y)Y^0u(y)^{-1}$ for an invertible element $u(y)$ of $\C(y)\otimes \End(V)$ such that $u(y)-1\in W_{-1}(\C(y) \otimes \End(V))$ and  $\delta_W(yN_1+N_2)\in \C(y)\otimes \gr^W\End(V)$. Hence by Theorem \ref{delthm2}, Theorem \ref{thm1var} is reduced to the following Proposition \ref{prop1var}. 

 Let $\tau_1:\mathbb{G}_m \to \Aut(V)$ be the homomorphism of algebraic groups corresponding to the splitting $Y^1$ of $W^1$ and let $t(y)= \tau_1(y^{-1/2})$. Consider $$N_y=t(y)^{-1}(yN_1+N_2)t(y)= N_1+ t(y)^{-1}N_2t(y).$$ 
 We have $$\spl_W(N_y)= t(y)^{-1}\spl_W(yN_1+N_2)t(y)=u^{\natural}(y) Y^0u^{\natural}(y)^{-1}\quad \text{with} \;u^{\natural}(y)= t(y)^{-1}u(y)t(y),$$ $$\delta_W(N_y)= t(y)^{-1}\delta_W(yN_1+N_2)t(y).$$
 
 By Theorem \ref{delthm2}, for  $y\gg 0$, we can write
$$u^{\natural}(y)= 1+ \sum_{m=1}^{\infty} u^{\natural}_m y^{-m/2}, \quad \delta_W(N_y)=\sum_{m=0}^{\infty} \delta^{\natural}_my^{-m/2}$$
as convergent Taylor series. 

\end{para}

\begin{prop}\label{prop1var} $(1)$ The $Y^1$-weights $s$ of $u_m^{\natural}$ satisfy $|s|<m$. 

$(2)$ We have $\delta^{\natural}_0=\delta_W(N_1)$. For $m\geq 1$, the $Y^1$-weights $s$ of $\delta_m^{\natural}$ satisfy $-m \leq s+2\leq m$. 
\end{prop}

We prove Proposition \ref{prop1var} in Lemma \ref{jk3}--\ref{Pfn2}. 

\begin{lem}\label{jk3} $(1)$ For every $a\in \Z$, we have $[N_1^{[a]}, \hat N_2]=0$ and $[N_1^{[a]}, \hat N_2^+]=0$. Here $N_1^{[a]}$ is the component of $N_1$ of $Y^0$-weight $a$.

$(2)$ For every $r\geq 1$, let $M^{(-r)}$ be the component of $N_2$ of $Y^1$-weight $-r$.  Then $[M^{(-r)}, \hat N_2^+]=0$.

\end{lem}

\begin{pf} This follows from Propositions \ref{jk1} and \ref{jk2}.  \end{pf}

\begin{para}\label{r=-1} We remark that $M^{(-1)}=0$. 
  In fact, since $Y^1$ is the Deligne splitting of $W^1$ with respect to $W^2$ and $N_2$, the part of 
$Y^1$-weight $-1$ of $N_2$ is zero (\ref{sdwt} (3)). 

\end{para}

\begin{lem}\label{filk} Let $D$ be a finite dimensional representation of $\SL(2)$ over $\C$. For $k\geq 0$, let  $\fil_kD\subset D$ be the sum of all subrepresentations of $D$ which are isomorphic to $\Sym^b(\rho)$ for $b\leq k$, where $\rho$ is the irreducible two dimensional representation of $\SL(2)$. Let
$$N=\begin{pmatrix} 0&1\\0&0\end{pmatrix}, \quad N^+= \begin{pmatrix} 0&0\\1&0\end{pmatrix}$$
which act on $D$ by the associated Lie action.
We have$:$

$(1)$ $\fil_0D=\{x\in D\;|\; N(x)=N^+(x)=0\}$. 

$(2)$ Let $k\geq 0$. If $N^+(x)=0$ and if $x$ is of $\SL(2)$-weight $\leq k$, then $x\in \fil_kD$. Here the $\SL(2)$-weight is defined by the action of $\mathbb{G}_m$ via $\mathbb{G}_m\to\SL(2)\;;\;a\mapsto \begin{pmatrix} a^{-1}& 0\\0 & a\end{pmatrix}$, that is, the weight defined by the eigenvalue of the Lie action of $\begin{pmatrix} -1 & 0\\ 0 & 1 \end{pmatrix}\in \sl(2)$. 

$(3)$ $\fil_kD \otimes \fil_{k'}D'\subset \fil_{k+k'}(D\otimes D')$. 

$(4)$ If $x\in \fil_kD$, then the $\SL(2)$-weights of $x$ is in the interval $[-k,k]$.
\end{lem}

\begin{pf}
These are evident.
\end{pf}

\begin{para} In what follows, we have  finite dimensional representations of $\mathbb{G}_m \times \SL(2)^2$. We consider $\fil_k$ for the action of the second $\SL(2)$. 

We consider a monodromy system $(V, W, N_1, N_2)$, the associated Deligne system $(N_1, N_2)$, and the associated action of $\mathbb{G}_m \times \SL(2)^2$. 
Let
$$N_y= t(y)^{-1}(yN_1+N_2)t(y)= N_1+\hat N_2+ \sum_{r\geq 1} y^{-r/2}M^{(-r)},$$ where $M^{(-r)}$ is as in Lemma \ref{jk3}.
\end{para}

\begin{lem}
\label{l:fil_property}
$(1)$ $N_1\in \fil_0\End(V)$. 

$(2)$ $\hat N_2 \fil_kV \subset \fil_kV$ and $[\hat N_2, \fil_k\End(V)]\subset \fil_k \End(V)$.

$(3)$ $M^{(-r)}\in \fil_{r-2}\End(V)$ for $r \geq 2$. 
\end{lem}

\begin{pf} (1) is by Lemma \ref{jk3} (1) and Lemma \ref{filk} (1). (2) holds because each subrepresentation of the second $\SL(2)$ is stable under $\hat N_2$.

We prove (3). $M^{(-r)}$ is $Y^1$-weight $-r$ and $Y^2$-weight $-2$ and hence 
  the $\SL(2)$-weight of $M^{(-r)}$ for the action of the second $\SL(2)$ (that is, the $\tau_1^{-1}\tau_2$-weight,  or the $Y^2-Y^1$-weight) is $r-2$. Since  $[M^{(-r)}, \hat  N_2^+]=0$ (Lemma \ref{jk3} (2)), (3) follows from Lemma \ref{filk} (2).    \end{pf}

\begin{lem}
\label{l:tau_1weights} Let $D$ be either $V$ or $\End(V)$, let
 $v\in \fil_kD$, and assume that $v$ is in the part of $Y^2$-weight $h$. 
  Then 
the $Y^1$-weights $s$ of $v$ satisfy $-k\leq s-h\leq k$. For example, in the case $h=0$ (resp.\ $h=-2$), we have $-k\leq s\leq k$ (resp.\ $-2-k\leq s\leq -2+k)$. (We will use these cases.) 

\end{lem}

\begin{pf} $h-s$ is the $Y^2-Y^1$ weight and it is the $\SL(2)$-weight of the action of the second $\SL(2)$. Hence this follows from Lemma \ref{filk} (4).
\end{pf}

\begin{para}\label{cM} Let $\cR:=\R[[y^{-1/2}]]$.
For a representation $D$ of $\mathbb{G}_m \times \SL(2)^2$, let $\cM D$ be the $\cR$-submodule of $\cR \otimes_{\C} D$  generated by 
  $\fil_0$ and $y^{-k/2-1}\fil_k$ for $k\geq 1$.

\end{para}

\begin{lem}\label{filfil} 
The maps $\End(V) \times V \to V$ and $\End(V) \times \End(V) \to \End(V)$ send $\fil_k \times \fil_{k'}$ to $\fil_{k+k'}$ and send (after taking $\cR \otimes$) $\cM \times \cM$ to $\cM$. 

\end{lem}

\begin{pf} This follows from Lemma \ref{filk} (3).
\end{pf}

\begin{lem}\label{Nyfil}
For $D=V$ (resp.\ $D=\End(V)$), $N_y$ (resp.\ $[N_y,\,\cdot\,]$) sends $\cM D$  to $\cM D$. 
\end{lem}

\begin{pf} This follows from Lemmas \ref{l:fil_property} and \ref{filfil}.
\end{pf}

\begin{lem}
\label{l:cM} Let $w<0$ and let $D=V^{[w]}$ or $D=\End(V)^{[w]}$, where $(\,\cdot\,)^{[w]}$ denotes the component of weight $w$ for $Y^0$. 
Let $D\{k\}$ be the part of $D$ of $Y^2$-weight $k$. Let $\cN=N_y^{[0]}$ in the case $D=V^{[w]}$, and let $\cN=\Ad(N_y^{[0]})$ in the case $D=\End(V)^{[w]}$. Let $\mu=-w$. 
 Then $\cN^\mu: \cR \otimes D\{0\} \overset{\cong}\to \cR\otimes  D\{2w\}$ and this induces 
$\cM D \cap (\cR\otimes  D\{0\}) \overset{\cong}\to \cM  D \cap (\cR \otimes D\{2w\})$. 
\end{lem}

\begin{pf} The map $\cN^\mu: \cR \otimes D\{0\} \to \cR\otimes D\{2w\}$ is an isomorphism as is seen by taking modulo $y^{-1/2}$. By Lemma \ref{l:fil_property}, 
$\cN^\mu$ sends the left-hand-side of the second expected  isomorphism 
to the right-hand-side. We prove the surjectivity. 
 Let $A$ be the inverse of the isomorphism $(\hat N_1+\hat N_2)^\mu: D\{0\}\to D\{2w\}$, and let $B= \cN^{\mu}-(\hat N_1+\hat N_2)^{\mu}$.  Then the inverse of the isomorphism  $\cN^\mu: \cR \otimes D\{0\} \to \cR\otimes D\{2w\}$ is written
as $\sum_{j=0}^{\infty} (-1)^j (AB)^jA$ (because the right-hand-side is $(1+AB)^{-1}A= (A^{-1}+B)^{-1}$ and $A^{-1}+B= \cN^{\mu}$). Hence this inverse of the isomorphism
$\cN^\mu$ sends $\cM D\cap (\cR \otimes D\{2w\})$ to $\cM D\cap (\cR \otimes D\{0\})$.
\end{pf}

\begin{lem}\label{l:cM2}  Let $D$ and $\cN$ be as in Lemma $\ref{l:cM}$. Let $a$ be an element of $\cM D\cap (\cR \otimes D\{-2\})$. Then we have   $a=\cN(b)+p$ for a unique pair $(b,p)$ such that $b\in \cM D\cap (\cR \otimes D\{0\})$, $p\in  \cM D\cap (\cR\otimes  D\{-2\})$, and  such that  $p$ is primitive for $\cN$.

\end{lem}

\begin{pf} Since $\cN^{\mu-1}a \in \cM D \cap (\cR \otimes D\{2w\})$, Lemma \ref{l:cM} shows that there is a unique $b\in \cM D \cap (\cR \otimes D\{0\})$ such that $\cN^{\mu-1}a= \cN^{\mu}b$. Let $p=a- \cN b\in \cM D \cap (\cR \otimes D\{-2\})$. Then $\cN^{\mu-1}p=0$, that is, $p$ is primitive for $\cN$. The pair $(b,p)$ is unique because it must be given by these constructions. 
\end{pf}

\begin{para}\label{Y0y}
Let $Y^0_y$ be the Deligne splitting of $W$ given by $W^2$ and $N_y$. Let $(\,\cdot\,)^{[w]}$ (resp.\, $(\,\cdot\,)^{[w,y]}$) be the component of $Y^0$ (resp.\ $Y^0_y$)-weight $w$. 
Then for  $w\in \Z$, we have a commutative diagram
$$\begin{matrix} \End(V)^{[w]} &\overset{\cong}\to & \gr^W_w\End(V)\\
\downarrow &&\Vert\\
\End(V)^{[w,y]} &\overset{\cong}\to & \gr^W_w\End(V),\end{matrix}$$
where the left vertical arrow is $f\mapsto u^{\natural}(y) f u^{\natural}(y)^{-1}$.

\end{para}

\begin{para}\label{Pfn1} Let $p_y\in \cR \otimes\End(V)$ be the element corresponding to $\delta_W(N_y)\in \cR \otimes \gr^W\End(V)$ via the isomorphism $\End(V) \cong \gr^W\End(V)$ given by $Y^0$. Since $N_y\in \cR \otimes \End(V)$ corresponds to $\gr^W(N)+\delta_W(N)\in \gr^W\End(V)$ via $Y^0_y$, we have
$$N_y=  u^{\natural}(y)(N_y^{[0]}+p_y)u^{\natural}(y)^{-1}$$
by \ref{Y0y}. 
We have $p_y\in \cR \otimes W_{-2}\End(V)$, 
$p_y$ is of $Y^2$-weight $-2$, and $p_y^{[w]}$ for $w\leq -2$ belongs to the primitive part in $\cR \otimes \End(V)^{[w]}$ for $\Ad(N_y^{[0]})$. These properties characterize $p_y$. 

For the proof of Proposition \ref{prop1var} (1) (resp.\ (2)), 
by the case $h=0$ (resp.\ $h=-2$) of Lemma \ref{l:tau_1weights} and by the fact that $u^{\natural}(y)$ (resp.\ $p_y$) is of $Y^2$-weight $0$ (resp.\ $-2$), it is sufficient to prove the following statement.  
\begin{equation}\tag{\ref{Pfn1}.1}
\text{Let $w<0$. Then $u^{\natural}(y)^{[w]}$ and $p_y^{[w]}$ belong to $\cM\End(V)$.}
\end{equation}

We prove this by using the induction on $-w>0$.

\end{para}

\begin{para}\label{Pfn2} We prove (\ref{Pfn1}.1).  By taking the part of $Y^0$-weight $w$ of $N_y u^{\natural}(y)= u^{\natural}(y) (N_y^{[0]}+p_y)$, 
 we have 
 $$ N_y^{[0]} u^{\natural}(y)^{[w]}+ \sum_{0\ge k>w} N_y^{[w-k]} u^{\natural}(y)^{[k]} = u^{\natural}(y)^{[w]} N_y^{[0]} + p_y^{[w]} + \sum_{0>k>w}  u^{\natural}(y)^{[k]}p_y^{[w-k]},$$
 that is, $$N_y^{[0]}u^{\natural}(y)^{[w]} - u^{\natural}(y)^{[w]}N_y^{[0]} - p_y^{[w]} =  (\sum_{0>k>w}  u^{\natural}(y)^{[k]}p_y^{[w-k]})-(\sum_{0\geq k>w} N_y^{[w-k]} u^{\natural}(y)^{[k]}).$$

Applying Lemma \ref{l:cM2} to $D=\End(V)^{[w]}$ by taking $(\sum_{0>k>w}  u^{\natural}(y)^{[k]}p_y^{[w-k]})-(\sum_{0\geq k>w} N_y^{[w-k]} u^{\natural}(y)^{[k]})$ as $a$, which is in $\cM\End(V)^{[w]}$ by induction (use $\fil_k\End(V) \cdot \fil_{\ell}\End(V)\subset \fil_{k+\ell}\End(V)$ by Lemma \ref{filk} (3)), we obtain that $b=u^{\natural}(y)^{[w]}$ and $-p=p_y^{[w]}$ are in $\cM\End(V)^{[w]}$.

This completes 
the proof of Proposition \ref{prop1var} and hence the proof of Theorem \ref{thm1var}. 
\end{para}

\begin{para} For a vector space $V$ with an increasing filtration $W$ and for a linear operator $N:V\to V$, we say that $(W, N)$ {\it splits} if there is a splitting of $W$ which is compatible with $N$.

\end{para}

\begin{thm}\label{mild1var}  Let $(N_1, N_2)$ be a monodromy system and assume that $(W,N_1)$ splits.

$(1)$  For $y\gg 0$, $\delta_W(yN_1+N_2)$ is a convergent Taylor series in $y^{-1}$.

$(2)$ Let $(\tau_j)_{j=0,1,2}:\mathbb{G}_m^3 \to \Aut(V)$ be homomorphism of algebraic groups such that $\tau_0$ splits $W$, $\tau_1$ splits the relative monodromy filtration of $N_1$ with respect to $W$, $N_1$ is of weight $-2$ for $\tau_1$, and $\tau_2$ is given by $Y$ of this monodromy system. Let $\tau^{\star}_1:=\tau_1\tau_0^{-1}$,   $t^{\star}(y):= \tau_1^{\star}(y^{-1/2})$, and $\delta_W^{\star}(yN_1+N_2):= t^{\star}(y)^{-1}\delta_W(yN_1+N_2)t^{\star}(y)$. Then for $y\gg 0$, $\delta^{\star}_W(yN_1+N_2)$ is a convergent Taylor series in $y^{-1/2}$. 
\end{thm} 

The proof is given in \ref{pfmild1var} below.

\begin{lem}\label{l:spl} Let  $(W,N,Y)$ be a Deligne system. Assume  that $(W, N)$ splits. Then $N$ is compatible with the splitting $Y^0$ of $W$ associated to this Deligne system and we have $\delta_W(N)=0$.

\end{lem}

\begin{pf} Fix a splitting $(Y')^0$  
of $W$ which is compatible with $N$. For each $w\in \Z$, choose a splitting of the relative monodromy filtration of $N$ on $\gr^W_w$ for which $N$ on $\gr^W_w$ is of weight $-2$. Via $(Y')^0$, these splittings give a splitting $(Y')^1$ of the relative monodromy filtration of $N$ with respect to $W$. Then $(Y')^0$ 
coincides with the Deligne splitting of $W$ associated to the Deligne system $(W, N, (Y')^1)$ and the $\delta_W(N)\in \gr^W\End(V)$ given by $(Y')^1$ is zero. Since $\delta_W(N)$ given by $Y^1$ coincides with that given by $(Y')^1$ (\ref{spldel}), it is zero and $Y^0$ is compatible with $N$.\end{pf}

\begin{para}\label{pfmild1var} We prove Theorem \ref{mild1var}. 

(1) follows from Theorem \ref{thm1var} (2) because  $\delta_W(N_1)=0$ by Lemma \ref{l:spl}. 

We prove (2). Since $N_1$ is of weight $-2$ for $\tau^{\star}_1$, we have 
  $$t^{\star}(y)^{-1}(yN_1)t^{\star}(y)=N_1.$$ 
  Hence 
   $$\delta^{\star}_W(yN_1+N_2) = \delta_W(N_1+t^{\star}(y)^{-1}N_2t^{\star}(y)).$$
   Since $N_1$ is compatible with $\tau_0$ (Lemma \ref{l:spl}), $\tau^{\star}_1$ splits the monodromy filtration (not the relative monodromy filtration) of $N_1$. 
   Since  $[N_1, N_2]=0$, $N_2$ respects the monodromy filtration of $N_1$.
  Hence the  $\tau_1^{\star}$-weights of $N_2$ are $\leq 0$. Hence  $t^{\star}(y)^{-1}N_2t^{\star}(y)$ 
is a polynomial in $y^{-1/2}$.

By \ref{calD}, it is sufficient to prove that the limit of $t^{\star}(y)^{-1}(yN_1+N_2)t^{\star}(y)$ for $y\to \infty$ belongs to the space $\cD$ in \ref{calD}. In the case where $\tau_j$ are given by $Y^j$ of the Deligne system $(N_1, N_2)$, this limit is $N_1+\hat N_2$ and it belongs to $\cD$. In general, there is $g\in \Aut(V)$ which preserves $W^j$ for $j=0,1,2$ and $\tau_j$ is given by $gY^jg^{-1}$. Then the limit of $t^{\star}(y)^{-1}(yN_1+N_2)N_2t^{\star}(y)$ coincides with $g(N_1+\hat N_2)g^{-1}$ and the map $M \mapsto gMg^{-1}$ preserves $\cD$.

\end{para}

\section{SL(2)-orbit theorems on $\spl_W$}
\label{s:spl}

 In this section, we prove one of the two main theorems in this paper (Theorem \ref{splthm}).
 
\begin{thm}\label{splthm0} Let $(V, W, \sig, Y)$ be a monodromy system with $\sig\neq \{0\}$. 
Then the  map
$$\sig^{\circ}/\R_{>0} \to \spl(W)\;;\; N\mapsto \spl_W(N)$$ extends to a real analytic map $\sig_{[:]}\to \spl(W)$ in the narrower sense. 
\end{thm}

Note that $\spl_W(N)$ is defined because $(W,N,Y)$ is a Deligne system, which is by 
Corollary \ref{c:adm} (1).

This theorem is equivalent to (1) of the following Theorem \ref{splthm}.

\begin{thm}\label{splthm} Let $(V, W, \sig, Y)$ be a monodromy system. Let the situation be as in $\ref{Njk2}$. Let $W^j$ and the splitting  $Y^j$ of $W^j$ be as in $\ref{ofp}$ associated to $p$. 

$(1)$ We have
$$\spl_W\Bigl(\sum_{(j,k)\in J} y_{j,k}N_{j,k}\Bigr)= u(y)Y^0u(y)^{-1},$$
where $u(y)$ is a convergent Taylor series in $y_{j+1,1}/y_{j,1}$ ($1\leq j\leq n-1$) and  $y_{j,k}/y_{j,1}  -c_{j,k}$ ($(j,k)\in J_1$)  with coefficients in $\End(V)$ such that the constant term of $u(y)$ is $1\in \End(V)$,  the coefficients of $u(y)$ except the constant term belong to $W_{-1}\End(V)$, and all coefficients of $u(y)$ are of $Y^n$-weight $0$.

$(2)$ $u(y)$ is in fact a rational function in $y_{j,k}$ ($(j,k)\in J$). 

$(3)$ Write $u(y)= \sum_{m\in \N^{n-1}} u_m \prod_{j=1}^{n-1} (y_{j+1,1}/y_{j,1})^{m(j)}$, where $u_m$ is in the ring $\cR$ af convergent Taylor series in $y_{j,k}/y_{j,1}-c_{j,k}$  ($(j,k) \in J_1$). Then $u_m \in \cR \otimes_\R \bigcap_{j=1}^{n-1} W^j_{\max(m(j)-1,0)}$.

\end{thm}

\begin{para}\label{Hodge1} The above Theorem \ref{splthm0} (resp.\ \ref{splthm}) is an analogue of the following  (0) (resp.\ (1)--(3)) in Hodge theory proved in \cite{KNU08}.

Consider a nilpotent orbit $(H_\R, W, \sig, F)$ with polarizable $\gr^W$.

(0) The map  $(\sig \times \R_{\geq 0})^{\circ}/\R_{>0}\cong \sig^{\circ} \to \spl(W)\;;\;\sig^{\circ}\ni  N \mapsto \spl_W(\exp(iN)F)$ extends to a real analytic function on some open neighborhood of $\sig_{[:]}$ in $(\sig\times \R_{\geq 0})_{[:]}$. 

 Let the situation be as in \ref{Njk2} and assume further  that $y_{n,1}\to \infty$.  (This tells that we assume $\sum_{(j,k)\in J} y_{j,k}N_{j,k} \in \sig^{\circ}\subset (\sig \times \R_{\geq 0})_{[:]}$  ($y_{j,k}\in \R_{>0}$) is sufficiently near to $p$ in $(\sig\times \R_{\geq 0})_{[:]}$.) 
  Consider the monodromy system $(H_\C, W_\C, \sig, Y)$ associated to this nilpotent orbit (\ref{Hratio}), and let $W^j$ and the splitting $Y^j$  of $W^j$ be as in \ref{ofp} associated to $p$.

(1) We have
$$\spl_W\Bigl(\exp\bigl(\sum_{(j,k)\in J} iy_{j,k}\bigr)F\Bigr) = u(y)Y^0u(y)^{-1},$$
where $u(y)$ is a convergent Taylor series in $y_{j+1,1}/y_{j,1}$ ($1\leq j\leq n$, $y_{n+1,1}$ denotes $1$) and  $y_{j,k}/y_{j,1}-c_{j,k}$ ($(j,k)\in J_1$)  with coefficients in $\End(V)$ with constant term $1$ such that all coefficients of $u(y)$ except the constant term belong to $W_{-1}\End(V)$.

 (2) $u(y)$ is in fact a rational function in $y_{j,k}$ ($(j,k)\in J$). 

(3) Write $u(y)= \sum_{m\in \N^n} u_m \prod_{j=1}^n (y_{j+1,1}/y_{j,1})^{m(j)}$, where $u_m$ are in the ring $\cR$ of  convergent Taylor series in $y_{j,k}/y_{j,1}-c_{j,k}$. Then $u_m\in \cR\otimes_\R \bigcap_{j=1}^n W^j_{\max(m(j)-1, 0)} \End(H_\R)$.

In fact, by the property of the real analytic structure of $(\sig\times \R_{\geq 0})_{[:]}$ around $\sig_{[:]}$ given in \ref{[:]'}, (0) and (1) are equivalent. (1)--(3) are essentially the contents of \cite{KNU08} Theorem 0.5 (1) and (2) (by the help of \cite{KNU08} 10.7).

\end{para}

\begin{para}\label{ton=2} We prove Theorem \ref{splthm}. 

(2) follows from the theory of Deligne systems over the field of rational functions in $y_{j,k}$. By this and by Theorem \ref{delthm2} (1) on $u^{\natural}(y)=t(y)^{-1}u(y)t(y)$, $u(y)$ is a Laurent series (obtained from a convergent Taylor series by multiplying $\prod_{j=1}^{n-1} (y_{j,1}/y_{j+1,1})^{a(j)}$ for some integers $a(j)\geq 0$).  
By this and by \ref{sdwt} (1) and (2), it is sufficient to prove that $u(y)$ is a Taylor series and the coefficients $u_m$ have the stated property on the $W^j$-weights for $1\leq j\leq n-1$. 

The last things are reduced to the case of a monodromy system $(N_1, N_2)$ in Section \ref{s:1var} as follows (the argument continues to \ref{pl1}).

Fix $j$ such that $1\leq j\leq n-1$. For each $(i,k)\in J$, take $\alpha_{i,k}\in \R_{>0}$ satisfying the following conditions. (i) $\alpha_{n,1}=\alpha_{j,1}=1$. (ii) If $y\gg 0$, $\sum_{(i,k)\in J, i\leq j} y\alpha_{i,k}N_{i,k} + \sum_{(i,k)\in J, i>j} \alpha_{i,k}N_{i,k} \bmod \R_{>0}$ is sufficiently near to $p$. By Corollary \ref{c:adm} (1), we have  the monodromy system $(N_1, N_2)$, where $N_1= \sum_{(i,k)\in J, i\leq j} \alpha_{i,k}N_{i,k}$ and $N_2= \sum_{(i,k)\in J, i>j} \alpha_{i,k}N_{i,k}$.

\end{para}

\begin{para}\label{pl1}  Let $(Y')^0$, $(Y')^1$, $(Y')^2= Y^n$ be those associated to the Deligne system $(N_1, N_2)$. Then we have Theorem \ref{thm1var} (1) on $\spl_W(yN_1+N_2)$. For the reduction to this result, we have to compare $(Y^0, Y^j)$ with $((Y')^0, (Y')^1)$. 
We have $\spl_W(yN_1+N_2)= u(y, \alpha) Y^0 u(y,\alpha)^{-1}$, where $u(y, \alpha)$ is obtained from $u(y)$ by the evaluation $y_{i,k}= y\alpha_{i,k}$ for $i\leq j$ and $y_{i,k}= \alpha_{i,k}$ for $i>j$. On the other hand, by Theorem \ref{thm1var} (1), we have $\spl_W(yN_1+N_2)= u'(y) (Y')^0 u'(y)^{-1}$, where $u'(y)$ denotes $u(y)$ in Theorem \ref{thm1var} (1). There is an element $g$ of $\Aut(V)$ such that $(Y')^0= gY^0g^{-1}$ and $(Y')^1= gY^jg^{-1}$, $g-1\in W_{-1}\End(V)+W^j_{-1}\End(V)$.  
We have a unique decomposition $g=g_1g_0$, where $g_1$ and $g_0$ respect $W$, $g_1-1\in W_{-1}\End(V)$, and $g_0$ commutes with $Y^0$. These $g_1$ and $g_0$ respect $W^j$.  
 We have 
  $u(y, \alpha)=u'(y)g_1$. Since $g_1\in W^j_0\End(V)$ and $u'(y)$ is a Taylor series whose $m$-th coefficient belongs to $W^j_{\max(m(j)-1,0)}\End(V)$, we have that $u(y, \alpha)$ is a Taylor series and its $m$-th coefficient belongs to $W^j_{\max(m(j)-1,0)}\End(V)$.  
  This proves Theorem \ref{splthm}. 
    \end{para}

\section{SL(2)-orbit theorems on $\delta_W$}
\label{s:del}

 In this section, we prove the other main theorem (Theorem \ref{delthm0}) and its variant (Theorem \ref{delthm3}).

\begin{para}\label{pl6} Let $\sig$ be a sharp finitely generated cone.  Assume $\sig\neq \{0\}$. 
Let $\cO^{\nar}$ be the sheaf of narrower real analytic functions on $\sig_{[:]}$. We define the invertible $\cO^{\nar}$-module $\frak L$ on $\sig_{[:]}$ as follows. 

Let $\sig'=\sig \times \R_{\geq 0}$ and let $(\cO')^{\text{nar}}$ be the sheaf of narrower real analytic functions on $(\sig')_{[:]}$. Recall that $(\sig')_{[:],<1}\smallsetminus \sig_{[:]}$ is an $\R_{>0}$-torsor over $\sig_{[:]}$ (\ref{U'}). We have an $(\cO')^{\nar}$-line bundle $\tilde{\frak L}=y_{1,1}(\cO')^{\nar}$ on $(\sig')_{[:],<1}\smallsetminus \sig_{[:]}$  by using the open covering $\{U(\Psi')\}_{\Psi}$ of $(\sig')_{[:],<1}$ in \ref{U'}
 ($y_{1,1}$ exists only locally but this line bundle glues globally). This descends to an  $\cO^{\nar}$-line bundle  $\frak L$  on $\sig_{[:]}$. For an open set $V$ of $\sig_{[:]}$, a section of this line bundle on $V$ is a section of $\tilde{\frak L}$ on the inverse image $\tilde V$ of $V$ in $(\sig')_{[:],<1}\smallsetminus \sig_{[:]}$ of degree $1$. Here degree $d$ means that it is multiplied by $c^d$ for the action of $c\in \R_{>0}$ on $(\sig')_{[:],<1}\smallsetminus \sig_{[:]}$.

\end{para}

In the following theorem, (1)  is a reformulation of (0).

\begin{thm}\label{delthm0} Let $(V, W, \sig, Y)$ be a monodromy system with $\sig\neq \{0\}$. 

$(0)$ $\delta_W$ is a section of the line bundle $\frak L$ on $\sig_{[:]}$.

Let the situation be as in $\ref{Njk2}$. 

$(1)$  
$$y_{1,1}^{-1}\delta_W\bigl(\sum_{(j,k)\in J} y_{j,k}N_{j,k}\bigr)$$ is a convergent Taylor series in $y_{j+1,1}/y_{j,1}$ ($1\leq j \leq n-1$) and $y_{j,k}/y_{j,1}  -c_{j,k}  $ ($(j,k)\in J_1$) with coefficients in $\gr^W\End(V)$. 

$(2)$ The constant term of this Taylor series is $\delta_W(\sum_{k=1}^{r(1)} c_{1,k}N_{1,k}) =\delta_W(N_1)$. All coefficients of this Taylor series belong to $W_{-2}\gr^W\End(V)$ and are of $Y^n$-weight $-2$.

$(3)$ The Taylor series in $(1)$  is in fact  a rational function in $y_{j,k}$. 

$(4)$ Write this Taylor series  in $(1)$ as $\sum_{m\in\N^{n-1}} a_m\prod_{j=1}^{n-1} (y_{j+1,1}/y_{j,1})^{m(j)}$, where $a_m$ is in the ring $\cR$ of  Taylor series in $y_{j,k}/y_{j,1}-c_{j,k}$ ($(j,k)\in J_1$). Then $$a_m \in \cR \otimes_\R \bigcap_{j=1}^{n-1} \; W^j_{m(j)-2}\gr^W\End(V).$$

\end{thm}

\begin{pf} 
 By Theorem \ref{delthm2}, we have that $\delta_W(N_y)= t(y)^{-1}\delta_W(\sum_{j,k} y_{j,k}N_{j,k})t(y)$ is a convergent Taylor series. Hence (locally on the space of ratios) $\delta$ is a Laurent series. By \ref{sdwt} (3) and (4) and by the argument in \ref{ton=2}, we are reduced to Theorem \ref{thm1var}  (2). 
\end{pf}

\begin{para} Let the notation be as in Theorem \ref{delthm0}. Consider  $\tau_j$ associated to $p$, and let  
\begin{align*}
\tau_j^{\star}&:=\tau_j\tau_0^{-1}, \\ 
t^{\star}(y)  &:=\prod_{j=1}^{n-1} \tau^{\star}_j((y_{j+1,1}/y_{j,1})^{1/2})=t(y)\tau_0((y_{1,1}/y_{n,1})^{1/2}), \\ 
\delta_W^{\star}(N) &:= t^{\star}(y)^{-1}\delta_W(N) t^{\star}(y). 
\end{align*}
\end{para}

\begin{thm}\label{delthm3} Let $(V, W, \sig, Y)$ be a monodromy system. Let the situation be as in $\ref{Njk2}$. Assume  that $(W, N)$ split for all $N\in \sig_{n-1}$.
Then 
$$\delta_W\bigl(\sum_{(j,k)\in J} y_{j,k}N_{j,k}\bigr)\quad (\text{resp.}\;\;  \delta^{\star}_W\bigl(\sum_{(j,k)\in J} y_{j,k}N_{j,k}\bigr))$$
with $y_{n,1}=1$ is a convergent Taylor series in $y_{j+1,1}/y_{j,1}$ (resp.\ $(y_{j+1,1}/y_{j,1})^{1/2}$)  ($1\leq j\leq n-1$) and $y_{j,k}/y_{j,1}-c_{j,k}$  ($(j,k)\in J_1$) with coefficients in $\gr^W\End(V)$. 

\end{thm}

\begin{pf} This is reduced to Theorem \ref{delthm0} and Theorem \ref{mild1var} (1) by  the method of \ref{ton=2}. \end{pf}

\begin{para}\label{delH3} The Hodge analogue of Theorem \ref{delthm3} is as follows.

Let $(\sig, F)$ be a nilpotent orbit with polarizable $\gr^W$. Let the situation be as in \ref{Njk2}.   Assume $n\geq 1$ and assume that $(W, N)$ split for all $N\in \sig_{n-1}$.
Then 
$$\delta_W\bigl(\exp\bigl(\sum_{(j,k)\in J} iy_{j,k}N_{j,k}\bigr)F\bigr) 
\quad (\text{resp.}\;\;  
\delta^{\star}_W\bigl(\exp\bigl(\sum_{(j,k)\in J} iy_{j,k}N_{j,k}\bigr)F\bigr))$$
with $y_{n,1}=1$ is a convergent Taylor series in $y_{j+1,1}/y_{j,1}$ (resp.\ $(y_{j+1,1}/y_{j,1})^{1/2}$)  ($1\leq j\leq n-1$), $y_{j,k}/y_{j,1}-c_{j,k}$ $((j,k)\in J_1$). 

This Hodge analogue follows from the mild SL(2)-orbit theorem \cite{KNU4} Theorem 5.1.2.

   In fact,  \cite{CKS}  4.65 (ii) and \cite{KNU08} Proposition 10.8 tell that the invariants of the SL(2)-orbit theorem depend on the parameters real analytically if the relative monodromy filtrations are fixed. We can apply it 
 to  the monodromy operators  $y_{j,1}^{-1}\sum_{k=1}^{r(j)} y_{j,k}N_{j,k}$ with parameters  ($1\leq j \leq n-1$) and to the Hodge filtration  $\exp(\sum_{k=1}^{r(n)} iy_{n,k}N_{n,k})F$ with parameters.

  \end{para}

\begin{para}\label{delH0}

The description of $\delta_W$ in Hodge theory becomes more complicated than Theorem \ref{delthm0} and Theorem \ref{thm1var} (2)  in general. For a nilpotent orbit $(N,F)$ with polarizable $\gr^W$, $\delta_W(\exp(iyN)F)$ ($y\gg 0$) need not be of the form $a_{-1}y+a_0+\sum_{m=1}^{\infty} a_my^{-m}$  but becomes $P(y)+ \sum_{m=1}^{\infty} a_my^{-m}$, where $P(y)$ is a polynomial in $y$ whose degree can be greater than $1$. 
\end{para}

\section{Asymptotic behaviors of regulators and local height pairings}\label{s:reg}

In this section, we prove Theorem \ref{thmreg} and Theorem \ref{thmht} which describe the asymptotic behaviors of the regulator and the local height pairing in the algebraic geometry over a non-archimedean local field. These theorems are obtained from Theorems \ref{thm1var} and \ref{mild1var} on monodromy systems.

\begin{para}\label{reg2} We review the weight-monodromy conjecture. 

Let $K$ be a complete discrete valuation field  with residue field $k$. Let $\ell$ be a prime number which is different from the characteristic of $k$.

Let $X$ be a proper smooth scheme over $K$. Let $m\geq 0$.

 Consider $N:H^m(X)_{\ell} \to H^m(X)_{\ell}(-1)$ (\ref{Ndef}). 
 
We consider the weight-monodromy conjecture assuming either one of the following (i) and (ii). 

\smallskip

(i) $k$ is a finite field.

(ii) $X$ has a proper model $\frak X$ over the valuation ring $O_K$ of $K$ with a strict semi-stable reduction.

\smallskip

In the case (i), the weight-monodromy conjecture for $H^m(X)_{\ell}$ is that if $M$ denotes the monodromy filtration of $N$ on $H^m(X)_{\ell}$, the weight of the Frobenius action on $\gr^M_w$ is $w+m$ for every $w$.

Assume (ii). 
Let $(Y_\la)_{\la\in \La}$ be the irreducible components of $\frak X \otimes_{O_K} k$. Then we have the Rapoport--Zink spectral sequence (\cite{RZ})
$$E^{-i,m+i}_1=
\bigoplus_{j \ge \max(0,-i)}\bigoplus_{\substack{E \subset \Lambda \\ 
\sharp E=i+2j+1}}H^{m-i-2j}(Y_E \otimes_k \overline k, \bQ_{\ell}(-i-j))
\Rightarrow E^m_{\infty}= H^m(X)_{\ell},$$
where $Y_E = \bigcap_{\lambda \in E} Y_{\lambda}$. 
The Rapoport--Zink spectral sequence degenerates at $E_2$ (\cite{Na}).
The weight-monodromy conjecture is that the monodromy filtration of $N$ on $H^m(X)_{\ell}$ coincides with the filtration as an end term of this spectral sequence. 

In the case where $k$ is finite and $X$ has $\frak X$, the two formulations of the weight-monodromy conjecture are equivalent.

\medskip

\noindent 
Remark 1. In the case $\text{char}(k)=\text{char}(K)$, the weight-monodromy conjecture is known to be true. 
In fact, 
Ito (\cite{Ito}) gives a proof by developing the previous work of Deligne \cite{De} Th\'eor\`eme 1.8.4.  
  In the case $\text{char}(k)=\text{char}(K)=0$, there is a Hodge-theoretic proof by 
using works of  Deligne--Saito--Steenbrink. 

\medskip
\noindent 
Remark 2. The weight-monodromy conjecture  is true in the case 
$\dim(X) \leq 2$ by \cite{RZ} and in many cases by \cite{SP}. 
\end{para}

\begin{para}\label{reg1}  We consider the regulator map in the non-archimedean geometry which we will define  assuming the weight-monodromy conjecture.

Let $F$ be a field, let $X$ be a  proper  smooth scheme over $F$, and let  $Z\in \gr^rK_n(X)$ for some integers $r,n>0$. Here $\gr^r$ is taken for the $\gamma$ filtration. Let  $\ell$ be a prime number which is different from the characteristic of $F$.

Then we have an $\ell$-adic  representation $H_Z$ of $G_F:=\Gal(\bar F/F)$ with an exact sequence 
$$0\to H^m(X)_{\ell}(r) \to H_Z \to \Q_{\ell} \to 0\quad \text{with}\;\; m=2r-n-1.$$ Here $H^k(X)_{\ell}$ for $k\in \Z$ is the $\ell$-adic \'etale cohomology $H^k_{\et}(X\otimes_F \bar F, \Q_{\ell})$. In fact, we have the Chern class homomorphism $\gr^rK_n(X) \to H^{2r-n}_{\et}(X, \Q_{\ell}(r))$. We have the exact sequence 
$0\to H^1(G_F, H^m(X)_{\ell}(r)) \to H^{2r-n}_{\et}(X, \Q_{\ell}(r)) \to H^0(G_F, H^{2r-n}(X)_{\ell}(r))$ and the composition  
$\gr^rK_n(X) \to H^{2r-n}_{\et}(X, \Q_{\ell}(r))\to H^0(G_F, H^{2r-n}(X)_{\ell}(r))$ is zero as is seen by the reduction to the case where $F$ is a finitely generated field over the prime field in which $H^0(G_F, H^{2r-n}(X)_{\ell}(r))=0$ by Weil conjecture proved by Deligne. We have 
 $H^1(G_F, H^m(X)_{\ell}(r))=\Ext^1_{G_F}(\Q_{\ell}, H^m(X)_{\ell}(r))$. Let $H_Z$ be the extension corresponding to the image of $Z$ in this ext group.

Let $W$ be the increasing filtration on $H_Z$ defined as follows. $W_w$ is $0$ if $w<m-2r=-n-1$, is $H^m(X)_{\ell}(r)$ if $m-2r\leq w\leq -1$, and is $H_Z$ if $w\geq 0$.

Assume now that $F$ is $K$ in \ref{reg2} and assume that $\ell$ is different from the characteristic of the residue field $k$ of $K$. We have the nilpotent operator $N: H_Z\to H_Z(-1)$. 
 Assume that we are either in the case (i) or in the case (ii) in \ref{reg2} and assume the weight-monodromy conjecture for $H^m(X)_{\ell}$.
The regulator of $Z$ is defined as 
$$\text{reg}(Z) = \delta_W(N)\in H^m(X)_{\ell}(r-1).$$

Here in the case (i), $\delta_W(N)$ is defined because Frobenius weight 
filtration $W^1$ is the relative  monodromy filtration of $N$ with respect to $W$ and the Frobenius weight splitting (using a splitting of $G_K\to G_k$) is a splitting of  $W^1$ which is compatible with $W$ and for which $N$ is of weight $-2$. 
In the case (ii), $\delta_W(N)$ is defined for  the following reason. 
  By spreading out, 
we are reduced to the case where $k$ is finitely generated over the prime field, and then by Rapoport--Zink spectral sequence 
 (using a splitting of the projection from the tame quotient of $G_K$ to $G_k$), we see that $H_Z$ is a direct sum of sub-representations of pure Frobenius weights. (Frobenius weight is the weight of the Frobenius of closed points of a suitable integral scheme $S$ of finitely type over $\Z$ whose function field is $k$.) The Frobenius weight filtration $W^1$ is the relative monodromy filtration of $N$ on $H_Z$, and we have the above direct sum decomposition is a splitting of $W^1$ which is compatible with $W$ and for which $N$ is of weight $-2$.

\end{para}

\begin{para}\label{htdef} We consider the local height pairing in the non-archimedean geometry which we will define  assuming the weight-monodromy conjecture.

Let $F$ be a field and  let $X$ be a  proper  smooth scheme over $F$. Let  $\ell$ be a prime number which is different from the characteristic of $F$. Let $A,B$ on $X$ be algebraic cycles of codimensions $r,s$ with $r+s=d+1$ ($d=\dim(X)$),  the supports of $A$ and $B$ are disjoint, and $A, B$ are homologically  equivalent to $0$ for the $\ell$-adic \'etale cohomology $H^{2r}(X)_{\ell}(r)$ and $H^{2s}(X)_{\ell}(s)$, respectively.
 Then for this pair $Z=(A,B)$, the Galois representation $H_Z$ is defined  (Beilinson \cite{Be}, Bloch \cite{Bl}). This $H_Z$ has $W$ whose $\gr^W_w$ are $\Q_{\ell}$, $H^{2r-1}(X)_\ell(r)$, $\Q_{\ell}(1)$ for $w=0,-1,-2$ and are $0$ otherwise.  
  The exact sequence $0\to H^{2r-1}(X)_{\ell}(r) \to W_0/W_{-2}\to \Q_{\ell}\to 0$ is given by the class of $A$, and the exact sequence $0\to \Q_{\ell}(1) \to W_{-1}\to H^{2r-1}(X)_{\ell}(r)\to 0$ is given by the class of $B$.

Assume now that $F$ is $K$ in \ref{reg2} and assume that $\ell$ is different from the characteristic of the residue field $k$ of $K$.
We have the nilpotent operator $N: H_Z\to H_Z(-1)$.

 Assume that we are either in the case (i) or in the case (ii) in \ref{reg2} and assume the weight-monodromy conjecture for $H^{2r-1}(X)_{\ell}$.
The local height pairing $\langle A, B\rangle$ is defined as 
$$\langle A, B\rangle = \delta_W(N)\in \Q_{\ell}.$$ 
Here $\delta_W(N)$ belongs to $\Q_{\ell}$ because $\Hom(\gr^W_0H_Z, \gr^W_{-2}H_\Z(-1))= \Q_{\ell}$. This $\delta_W(N)$ is defined for the reasons described in \ref{reg1} by using the Frobenius weight filtration and its splitting of a Galois representation of a finitely generated field over the prime field.

\end{para}

\begin{para}\label{antop} In the following, for a complete discrete valuation field $K$ and for a scheme $S$ over $K$ of finite type, we consider the topology of the set $|S|$ of all closed points of $S$ called the analytic topology. 
(We call so to distinguish it from the algebraic one, Zariski topology.)

Assume first that $S$ is affine. Then it is 
the weakest topology such that the map $|S|\to \R\;;\;p\mapsto |f(p)|$ is continuous for every $f\in \cO(S)$. Here we fix a multiplicative valuation $|\,\cdot\,|$ of $K$  and extend it to finite extensions of $K$ uniquely. 
If $\Omega$ denotes the algebraic closure of $K$ and we endow $S(\Omega)$  with the weakest topology such that $S(\Omega)\to \Omega\;;\;p\mapsto f(p)$ is continuous for every $f\in \cO(S)$, it is the quotient topology of the topology of $S(\Omega)$ for the canonical surjection $S(\Omega)\to |S|$. 

For a general $S$, the analytic  topology of $|S|$ is the unique topology such that for every affine open set $U$ of $S$, $|U|$ is open in $|S|$ and equipped with the above topology.

\end{para}

\begin{para}\label{reg3}  In the following Theorems \ref{thmreg} and \ref{thmht}, we consider the following situation.

 Let $K$ be a complete discrete valuation field with finite residue field $k$.

 Let $S$ be a smooth connected curve over $K$ (which need not be proper) and let $X\to S$ be a proper morphism of schemes. For  $p\in |S|$ (\ref{antop}), let $K(p)$ be the residue field of $p$, which is a finite extension of $K$, and let $X(p)= X \times_S p$. 
 
Let $p_0\in |S|$. Assume that $X\smallsetminus X(p_0) \to S\smallsetminus \{p_0\}$ is smooth, and $X\to S$ is of strict semi-stable reduction at $p_0$.

   In Theorem \ref{thmreg}, we assume that we are given
    $Z\in \gr^rK_n(X\smallsetminus X(p_0))$ ($r,n>0$).

   In Theorem \ref{thmht}, let $d$ be the relative dimension of $X$ over $S$, let $1\leq r\leq d$, let $A$ and $B$ be algebraic cycles on $X$ which are flat over $S$ such that $A$ is of codimension $r$ and $B$ is of codimension $d+1-r$, and assume that the supports of $A$ and $B$ are disjoint
and that $A$ and $B$ are  fiberwise homologically equivalent to $0$.  
      
   In Theorem \ref{thmreg} (resp.\ \ref{thmht}), we assume that 
    the following 1 for $m=2r-n-1$ (resp.\ $m=2r-1$) and the following 2 are satisfied. 

   \medskip

1. $H^m(X(p)/p)_{\ell}:=H^m_{\et}(X\otimes_{K(p)} \bar{K(p)}, \Q_{\ell})$ 
satisfies the weight-monodromy conjecture of the case (i) in \ref{reg2} if $p\in |S|\smallsetminus \{p_0\}$ is sufficiently near to $p_0$ for the analytic topology of $|S|$ (\ref{antop}). 

2. Let $(Y_{\la})_{\la\in \La}$ be the irreducible components of $X(p_0)$. For every non-empty subset $E$ of $\La$ and for every $i$, $H^i(Y_E)_{\ell}$  satisfies ($Y_E=\bigcap_{\la \in E} Y_{\la}$ here is regarded as a scheme over the local field $K(p_0)$) the weight-monodromy conjecture of the case (i) in \ref{reg2}. 

\medskip

In Theorem \ref{thmreg}, for $p\in |S|$, $p\neq p_0$,  let $Z(p)\in \gr^rK_n(X(p))$ be the pullback of $Z$, and consider the regulator $\mathrm{reg}(Z)(p):=\delta_W(N(p))$, where $N(p)$ is the monodromy operator of the local field $K(p)$ with respect to $K$ (\ref{Ndef2}) on $H_{Z(p)}$ (\ref{reg1}). That is, $\mathrm{reg}(Z)(p)=e(K(p)/K)^{-1}\mathrm{reg}(Z(p))$, where the last term is the regulator in \ref{reg1} of $Z(p)$ over the local field $K(p)$ and $e(K(p)/K)$ is the ramification index.
  Let $L$ be the field of fractions of the completion of the discrete valuation ring  $\cO_{S,p_0}$. We denote the regulator of the pullback to $X_L:=X \times _S\Spec(L)$ of  $Z$ by $\text{reg}_{\text{geo}}(Z)$
    (note that $L$ is a local field of  equal characteristic and hence the weight-monodromy conjecture \ref{reg2} (ii) is true). 

   In Theorem \ref{thmht},  for $p\in |S|$, $p\neq p_0$, let $Z(p)=(A(p), B(p))$ be the pullback of $Z=(A,B)$ to $X(p)$, and consider the local height pairing  $\langle A, B\rangle(p):=\delta_W(N(p))$, where $N(p)$ is the monodromy operator of the local field $K(p)$ with respect to $K$ (\ref{Ndef2}) on $H_{Z(p)}$ (\ref{htdef}). That is, $\langle A, B\rangle(p) =e(K(p)/K)^{-1}\langle A(p), B(p)\rangle$.  
   We denote the local height pairing of the pullback of $Z=(A, B)$ to  $X_L=X \times _S\Spec(L)$ 
by $\langle A,B\rangle_{\text{geo}}$. 

We will identify  $H^m(X(p)/p)_{\ell}$ for all $p\in |S|\smallsetminus \{p_0\}$ and $H^m(X_L/L)_{\ell}:=H^m_{\et}(X_L\otimes_L \bar L, \Q_{\ell})$ by identifying them with the stalk of an $\ell$-adic sheaf $H^m(X/S)_{\ell}$ on the log point $0$ in \ref{pf2}. We identify $H^m(X/S)_{\ell}$ with its stalk. Furthermore, we 
fix an isomorphism $\bar \Q_{\ell}\cong \C$, and we denote $H^m(X/S)_{\ell}\otimes_{\Q_{\ell}} \C$ by $H^m(X/S)_{\ell,\C}$.
Our regulators now have values in $H^m(X/S)_{\ell,\C}(r-1)=H^m(X/S)_{\ell,\C}$ (the last equality comes from 
$\Q_{\ell}(1)\otimes_{\Q_{\ell}} \C= (\Q_{\ell} \otimes_{\Q} \Q\cdot (2\pi i)) \otimes_{\Q_{\ell}} \C = \C$)
 and our local height pairings have values in $\Q_{\ell}(1) \otimes_{\Q_{\ell}} \C=\C$.

In Theorems \ref{thmreg} and \ref{thmht}, we fix a local parameter $t$ on $S$ at $p_0$, and we define $$y(p):= v(t(p)),$$ where $v$ is the additive valuation $K(p)^\times \to \Q$  of $K(p)$ which sends a prime element of $K$ to $1$. 
\end{para}

\begin{thm}\label{thmreg} 
Let the situation be as in $\ref{reg3}$.

$(1)$ We have 
 $$\mathrm{reg}(Z)(p)= \mathrm{reg}_{\mathrm{geo}}(Z)\cdot y(p)+\sum_{m=0}^{\infty} a_my(p)^{-m}$$ 
for $p\in |S|  \smallsetminus \{p_0\}$ which is near  to $p_0$ for the analytic topology ($\ref{antop}$), 
where $a_m\in H^m(X/S)_{\ell,\C}$ and $\sum_{m=0}^{\infty}$ is a convergent series (convergent for the topology of $\C$, not the $\ell$-adic topology).

$(2)$ If $Z$ comes from $\gr^rK_n(X)$, then $\mathrm{reg}_{\mathrm{geo}}(Z)=0$.

\end{thm}

\begin{thm}\label{thmht} Let the situation be as in $\ref{reg3}$. We have 
$$\langle A, B\rangle(p)= \langle A, B\rangle_{\mathrm{geo}}\cdot y(p)+\sum_{m=0}^{\infty} a_my(p)^{-m}$$ 
for $p\in |S|  \smallsetminus \{p_0\}$ which is near  to $p_0$ for the analytic topology ($\ref{antop}$), 
where $a_m\in \C$ and $\sum_{m=0}^{\infty}$ is a convergent series.

\end{thm}

\begin{para}\label{pf(2)}
We first prove Theorem \ref{thmreg} (2).  

Let $L_{\text{ur}}\subset \bar L$ be the maximal unramified extension of $L$ 
and let $O_{L_{\text{ur}}}$ be its valuation ring. We have a commutative diagram
$$\begin{matrix} \gr^rK_n(X) && \to && H^{2r-n}(X\times_S \Spec(O_{L_{\text{ur}}}), \Q_{\ell})(r)\\
\downarrow&&&&\downarrow \\
\gr^rK_n(X\smallsetminus X(p_0)) &\to &H^1(G_{L_{\text{ur}}}, H^m(X_L/L)_{\ell}(r)) &\to & H^{2r-n}(X\times_S \Spec(L_{\mathrm{ur}}), \Q_{\ell})(r) \\
&&\Vert && \downarrow \\
&& H^m(X_L/L)_{\ell}(r-1)/NH^m(X_L/L)_{\ell}(r)&\overset{\subset}\to & H_{2d-m}(X(p_0))_{\ell}(r-1-d).
\end{matrix}$$
Here: The upper horizontal arrow is the Chern class map. At the lower right corner, $d=\dim(X(p_0))$ and $H_{2d-m}(X(p_0))_{\ell}$ is the dual of $H^{2d-m}(X(p_0)) _\ell$.  The column on the right-hand-side is by the Poincar\'e duality theory for \'etale cohomology and is exact. $N$ at the bottom is the monodromy operator of the local field $L$. The middle vertical identification is due to the fact that the action of $G_{L_{\text{ur}}}$ on $H^m(X_L/L)_{\ell}$ is unipotent because $X\times_S \Spec(O_{L_{\text{ur}}})$ is of semi-stable reduction over $O_{L_{\text{ur}}}$, and it sends the class of the extension $0\to H^m(X_L/L)_{\ell}(r) \to H \to \Q_{\ell}\to 0$ to the class of $N\tilde 1$, where $\tilde 1$ is a lifting of $1\in \Q_{\ell}$ to $H$. 
The horizontal arrow in the bottom is injective because it is the dual of the map
$H^{2d-m}(X(p_0))_{\ell}(d-r+1)\to H^{2d-m}(X_L/L)_{\ell}^{N=0}(d-r+1)$ which is surjective by the local invariant cycle theorem (\cite{De} Theorem 3.6.1).

By the diagram, the composition $\gr^rK_n(X) \to H_{2d-m}(X(p_0))_{\ell}(r-1-d)$ is the zero map. Hence by the injectivity of the horizontal arrow in the bottom, the map 
$\gr^rK_n(X) \to H^m(X_L/L)_{\ell}(r-1)/NH^m(X_L/L)_{\ell}(r)$ is the zero map. Hence if $Z$ comes from $\gr^rK_n(X)$, $(W, N)$ of $H_Z$ splits and hence $\text{reg}_{\text{geo}}(Z)=0$.

\end{para}

We start the proofs of Theorem \ref{thmreg} (1) and Theorem \ref{thmht}. We prove them together.

\begin{para}\label{pf00} 
We consider the relative versions for $X\smallsetminus X(p_0)\to S\smallsetminus \{p_0\}$ of the constructions in \ref{reg1} and 
\ref{htdef}.

 Let $H^m(X/S)_{\ell}= R^mf_*\Q_{\ell}$, where $f: X\smallsetminus X(p_0) \to S\smallsetminus \{p_0\}$.
We have $H_Z$ on the \'etale site of $S\smallsetminus \{p_0\}$ with $W$.  $W_0H_Z=H_Z$. All non-zero $\gr^W$ are $\gr^W_0H_Z=\Q_{\ell}$ and $\gr^W_{-n-1}H_Z= H^m(X/S)_{\ell}(r)$ in the case of Theorem \ref{thmreg}, and $\gr^W_0H_Z=\Q_{\ell}$, $\gr^W_{-1}H_Z=H^m(X/S)_{\ell}(r)$, $\gr^W_{-2}H_Z=\Q_{\ell}(1)$  in the case of Theorem \ref{thmht}. This $H_Z$ is obtained as follows.

In the case of Theorem \ref{thmreg}, we use the Chern class map 
$$\gr^rK_n(X\smallsetminus X(p_0)) \to H^{2r-n}(X\smallsetminus X(p_0), \Q_{\ell}(r)) \supset H^1(S\smallsetminus \{p_0\}, H^m(X/S)_{\ell}(r))$$
which induces 
$$\gr^rK_n(X\smallsetminus X(p_0)) \to H^1(S\smallsetminus \{p_0\}, H^m(X/S)_{\ell}(r))=\text{Ext}^1(\Q_{\ell}, H^m(X/S)_{\ell}(r)).$$
The case of Theorem \ref{thmht} is obtained by using the construction of Beilinson and Bloch (\ref{htdef}) considering the relative version for $X\smallsetminus X(p_0)\to S\smallsetminus \{p_0\}$. 

\end{para}

\begin{para}\label{pf1} 
There are a finite extension $K'$ of $K$, a smooth connected curve $S'$ over $K'$ with a $K'$-rational point $p'_0$, and a flat $K$-morphism $S'\to S$ which sends $p'_0$ to $p_0$ and $S'\smallsetminus \{p_0'\}$ to $S\smallsetminus \{p_0\}$, such that the image of $\pi_1(S'\smallsetminus \{p'_0\})$ in the automorphism group of the stalk of $H_Z$ is a pro-$\ell$ group. Let $O_{K'}$ be the valuation ring of $K'$ and let $k'$ be the residue field of $K'$. 
Let $t'$ be a prime element of the discrete valuation ring $\cO_{S', p'_0}$  and let $\pi'$ be a prime element of $K'$.

Let $L'$ be the field of fractions of the completion of $\cO_{S',p_0'}$ 
and let $O_{L'}$ be the valuation ring of $L'$.
 Then $O_{L'}= K'[[t']]\supset \{\text{convergent Taylor series}\}\supset \cO_{S',p_0'}$, $\{\text{convergent Taylor series}\}= \varinjlim_h O_{K'}[[t'/(\pi')^h]][1/\pi']$. Replacing $t'$ by $t'/(\pi')^h$ for $h\gg 0$, we have that $\Spec(O_{L'})\to S'$ factors as $\Spec(O_{L'})\to \Spec(O_{K'}[[t']][1/\pi'])\to S'$. 

If $p\in |S|$ is near to $p_0$ for the analytic topology, we have $p'\in |S'|$ lying over $p$ which is near to $p_0'$ for the analytic topology, and since $|t'(p')|<1$, $p'\to S'$ factors through $\Spec(O_{K'}[[t']][1/\pi'])$.

\end{para} 

\begin{para}\label{pf2}

Consider the pullback of $H_Z$ on $\Spec(O_{K'}[[t']][1/t', 1/\pi'])$.
 By the theory of Grothendieck--Murre on the tame fundamental group \cite{GT}, 
 the action of the kernel $I$ of $\pi_1(\Spec(O_{K'}[[t']][1/t', 1/\pi']))\to \Gal(\bar {k'}/k')$ on the stalk of $H_Z$ factors through the surjection $I\to \Z_{\ell}(1) \times \Z_{\ell}(1)$ given by the $\ell^n$-th roots of $t'$ and $\pi'$ for $n\geq 0$.

Hence  $H_Z$ on $\Spec(O_{K'}[[t']][1/t', 1/\pi'])$
comes from  a smooth $\Q_{\ell}$-sheaf (which we still denote by $H_Z$) on the log \'etale site of $\Spec(O_{K'}[[t']])$ whose log structure is given by $t'$ and $\pi'$.
Then $H_Z$ induces  a smooth $\Q_{\ell}$-sheaf  (which we still denote by $H_Z$) on the log \'etale site of the log point $0:=\Spec(k')$ whose log structure is given by $t'$ and $\pi'$. The subquotient $H^m(X/S)_{\ell}(r)$ of $H_Z$ also induces a smooth $\Q_{\ell}$-sheaf  (which we still denote by $H^m(X/S)_{\ell}(r)$) on the log \'etale site of $0$.

By the structure of $\pi_1^{\log}(0)$ (by the inner automorphism, the geometric Frobenius in $\Gal(\bar {k'}/k')$ acts on  the above  $\Z_{\ell}(1) \times \Z_{\ell}(1)$ as the multiplication by $q^{-1}$, where $q$ is the order of $k'$), the action of $\Ker(\pi_1^{\log}(0)\to \Gal(\bar{k'}/k')$ on the stalk of $H_Z$ is unipotent.

Hence by Proposition \ref{propreg} below, Theorem \ref{thmreg} (1) and Theorem \ref{thmht} are  reduced to Theorem \ref{thm1var} (2).

Let $(N_1, N_2)$ be the base of the monodromy cone of $0$ which is dual to $(t', \pi')$. We have $N_1,N_2: H_Z\to H_Z(-1)$ (we identify $H_Z$ with its stalk). Let $W^2$ be the Frobenius weight filtration of $H_Z$ on $0$ and let $Y$ be the splitting of $W^2$ obtained by fixing  a splitting of the projection $\pi^{\log}(0) \to \Gal(\bar {k'}/k')$. Then $N_1$ and $N_2$ are of weight $-2$ for $Y$. 
\end{para}

\begin{prop}\label{propreg} Replace $t'$ by $t'/(\pi')^h$ for $h\gg 0$. Then this $H_Z$ on $0$ gives an object of $\cA_0$. That is, it gives a monodromy system $(N_1, N_2)$.

\end{prop}

To prove this, we  give lemmas.

\begin{lem}\label{pfreg1}  Replace $t'$ by $t'/(\pi')^h$ for $h\gg 0$. Then for $H_Z$ on $0$, we have that for every $y\in \R_{\geq 0}$, $W^2$ is the relative monodromy filtration of $yN_1+N_2$ with respect to $W$.

\end{lem}

This is reduced to 
\begin{lem}\label{pfreg2}  Replace $t'$ by $t'/(\pi')^h$ for $h\gg 0$. Then for $H^m(X/S)_{\ell}(r)$ on $0$, we have that for every $y\in \R_{\geq 0}$, $W^2$ on $H^m(X/S)_{\ell}(r)$ is the relative monodromy filtration of $yN_1+N_2$ with respect to $W$.

\end{lem}

\begin{pf} It is sufficient to prove that for the Frobenius weight 
filtration $\cW$ on $H^m(X/S)_{\ell}$, $(yN_1+N_2)^i$ induces an isomorphism $\gr^{\cW}_{m+i}\to \gr^{\cW}_{m-i}$ for every $i\geq 0$ after we replace $t'$ by $t'/(\pi')^h$ for $h\gg 0$, that is, after we replace $N_2$ by $hN_1+N_2$ for  $h\gg 0$. 
By the assumption 1 in \ref{reg3} and by \ref{Natp} with $n=1$,
$(yN_1+N_2)^i:  \gr^{\cW}_{m+i} \to  \gr^{\cW}_{m-i}$ for each $i\geq 0$ is an isomorphism for infinitely many $y$. Hence 
 the polynomial $\text{det}((yN_1+N_2)^i): \text{det}(\gr^{\cW}_{m+i}) \to \text{det}(\gr^{\cW}_{m-i})$ in 
  $y$ is non-zero. Hence  this polynomial does  not have zero for $y\gg 0$. Hence it does not have zero for $y\geq 0$ if we replace $N_2$ by $hN_1+N_2$ (and hence $yN_1+N_2$ by $(y+h)N_1+N_2$) for $h\gg 0$. 
\end{pf}

\begin{lem}\label{pfreg3} On $0$, let $W^1$ on $H_Z$ be the relative weight filtration of $N_1$ with respect to $W$. Then $W^2$ on $H_Z$ is the relative monodromy filtration of $N_2$ with respect to $W^1$.

\end{lem}

This is reduced to

\begin{lem}\label{pfreg4}  On $0$, let $W^1$ on $H^m(X/S)_{\ell}(r)$ be the relative weight filtration of $N_1$ with respect to $W$. Then $W^2$ on $H^m(X/S)_{\ell}(r)$ is the relative monodromy filtration of $N_2$ with respect to $W^1$.

\end{lem}

\begin{pf} This follows from the assumption 2 in \ref{reg3}. In fact,  2 tells that in the Rapoport--Zink spectral sequence, $E^{p,q}_1$ with the pure weight filtration of weight $q$ and with the action of $N_2$ and $Y^2$ is a Deligne system of one variable. Since the category of Deligne systems is an abelian category, $E^{p,q}_2=E^{p,q}_{\infty}$ of pure weight $q$ with $N_2$ and $Y^2$ is a Deligne system in one variable. 
 Since $W^1$ is the filtration as an end term, this tells that the Frobenius weight filtration $W^2$ is the relative monodromy filtration of $N_2$ with respect to $W^1$. 
\end{pf}

\begin{para}\label{pfreg5} By Lemmas \ref{pfreg1} and \ref{pfreg3}, we have Proposition \ref{propreg}. This completes the proofs of Theorem \ref{thmreg} and Theorem \ref{thmht}.

\end{para}

\begin{para}
The Hodge version of Theorem \ref{thmreg} (2) is given in \cite{KNU4} Proposition 7.2.3.
\end{para}

\begin{para}

 Theorem \ref{thmht}  is a non-archimedean analogue of  the archimedean theory of Pearlstein (\cite{Pe} Theorem 5.19) obtained by using his $\SL(2)$-orbit theorem for mixed Hodge structures.
\end{para}

\begin{para}

The asymptotic behavior of the local height pairing in the non-archimedean algebraic geometry can be formulated also by using the intersection theory on the model of $X$ over $O_K$ (not using $\delta_W$ of a monodromy system). The authors learned this from Spencer Bloch. This will be discussed in \cite{Ka2}.

In \cite{Ka2}, the asymptotic behavior of the global height pairing will be discussed by generalizing  the formula of Tate in \cite{Ta}, and also, the asymptotic behavior of heights of motives will be discussed. 

\end{para}

\section{On the expectation 6.16 of Part I}
\label{s:linalg}

  Here we prove an expectation given in Part I, 6.16, which provides the first nontrivial example giving an affirmative answer to Question 6.3 of Part I. 

\begin{prop}
\label{p:linalg}
  Let $X \to Y$ be the degenerate Tate elliptic curve over the standard log point over a finite field 
as in Part {\rm{I}}, $6.4$. 
  Let $H$ be a pure and simple object of $\cA_X$. 
  Assume that the action of the log fundamental group of a log geometric fiber on a stalk of $H$ is unipotent. 
Then $H$ is the pullback of an object of $\cA_Y$. 
\end{prop}

  This means that the equivalent conditions in Part {\rm{I}}, Proposition $6.11$ $(3)$ are always satisfied.
  Hence, by Part {\rm{I}}, $6.16$, the statement in Part {\rm{I}}, Question $6.3$ is true for this $X \to Y$: 

\begin{prop}
Let $f: X\to Y$ be as in Proposition $\ref{p:linalg}$. 
Then the higher direct image functors $R^mf_*$ 
for the log \'etale sites send objects of $\cA^s_X$ to $\cA_Y^s$. 
\end{prop}

Recall that $\cA_X^s$ denotes the category of objects of $\cA_X$ whose $\gr^W$ for the weight filtrations $W$ are semi-simple, and $\cA_Y^s$ are defined in the similar way.

  In the rest of this section, we prove Proposition \ref{p:linalg}.

\begin{para}
  First we recall the notation. 
  We denote $Y$ also by $y=\Spec(k)$.
  Let $p$ be the characteristic of $k$ and $q$ the order of $k$.
We consider the logarithmic  fundamental group $\pi_1^{\log}(X)$. 
  Fix a section $s: Y \to X$ of $X \to Y$ and let $x$ be the image of $s$. 
  Considering the stalk at $\bar x(\log)$ lying over $\bar y(\log)$, a
smooth $\bar \Q_{\ell}$-sheaf on the log \'etale site $X_{\loget}$ of $X$ is identified with a representation of $\pi_1^{\log}(X)$ over $\bar \Q_{\ell}$. 

We have canonical exact sequences 

\begin{sbpara}\label{pi1} \;\; 
$1 \to \pi_1^{\log}(X\times_Y \bar y(\log)) \to \pi_1^{\log}(X) \to \pi_1^{\log}(y)\to 1$, 
\end{sbpara}

\noindent and

\begin{sbpara}\label{pi1fiber} \;\; $1 \to \hat\bZ'(1) \to \pi_1^{\log}(X\times_Y \bar y(\log)) \to \hat \bZ \to 1$, 
\end{sbpara}

\noindent 
where $\hat\bZ'(1)=\prod_{\ell'\not=p}\bZ_{\ell'}(1)$.

\noindent 
(To see the latter, we use the surjection from $\pi_1^{\log}$ of the geometric generic fiber to $\pi_1^{\log}(X\times_Y \bar y(\log))$. Cf.\ \cite{Hoshi} Theorem 2.)

  Let $\gamma_0$ be a topological generator of $\pi_1^{\log}(\overline y)\subset \pi_1^{\log}(y)$. 
  Let $F \in \pi_1^{\log}(y)$ be a lift of the $q$-th power map of $\Gal(\bar k/k)$ (which is a quotient of $\pi_1^{\log}(y))$.
  Lift $\gamma_0$ and $F$ by the section $\pi_1^{\log}(y) \to \pi^{\log}_1(X)$ induced by $s$ and denote them by the same symbols. 
  Then we have $$
F\gamma_0F^{-1}=\gamma_0^q$$
in $\pi_1^{\log}(X)$.  
  Further, there are an element $\gamma_1$ of $\pi_1^{\log}(X\times_Y \bar y(\log))
\subset \pi_1^{\log}(X)$ whose image in $\hat \bZ$ is a topological generator and a topological generator $\gamma_2$ of 
$\hat \bZ'(1)\subset \pi_1^{\log}(X)$ 
such that 
$$
F\gamma_2F^{-1}= \gamma_2^q, \quad 
\gamma_0\gamma_1\gamma_0^{-1}= \gamma_1\gamma_2^a, \quad 
\gamma_0\gamma_2\gamma_0^{-1}=\gamma_2,$$
where $a$ is the image of the $q$-invariant of $X$ in $(M_Y/\cO_Y^{\times})_{\bar y}\cong \bN$. 
  Let $b \in \hat \bZ'$ such that $F\gamma_1F^{-1}=\gamma_1\gamma_2^b$.
  Note that in many cases (for example, when $a$ is coprime to $q-1$), we may assume $b=0$ by replacing $F$ and $\gamma_1$, but not always. 
\end{para}

\begin{para}
  Let $N_j$ be the endomorphism induced by $\log\gamma_j$ 
(resp.\ $a\log\gamma_j$) on the stalk of $H$ for $j=0,1$ (resp.\ $j=2)$. 
  Then we have $[N_0,N_2]=[N_1,N_2]=0$, $FN_0=qN_0F$, and 
$FN_2=qN_2F$. 
  Further, we have $[N_0, N_1]=N_2$, which is proved as follows.
  By $\gamma_0\gamma_1\gamma_0^{-1}= \gamma_1\gamma_2^a$, we have $\gamma_0 \gamma_1^n\gamma_0^{-1}= (\gamma_1\gamma_2^a)^n$ for every $n\geq 0$.
  Hence 
(below we denote the automorphism on the stalk of $H$ induced by an element of $\pi_1^{\log}$ by the same symbol) 
$$ \gamma_0 N_1\gamma_0^{-1} = \lim_{n\to \infty} \gamma_0 \frac{\gamma_1^{\ell^n}-1}{\ell^n}\gamma_0^{-1}= \lim_{n\to \infty} \frac{(\gamma_1\gamma_2^a)^{\ell^n}-1}{\ell^n}=\log(\gamma_1\gamma_2^a)= N_1+N_2.$$
  From this, we have $\gamma_0^n N_1 \gamma_0^{-n}= N_1+nN_2$ for every $n\geq 0$. 
  Therefore $$[N_0, N_1]= \lim_{n\to \infty} \frac{\gamma_0^{\ell^n}N_1\gamma_0^{-\ell^n}-N_1}{\ell^n} $$
$$=\lim_{n\to \infty} \frac{N_1+\ell^nN_2-N_1}{\ell^n} = N_2.$$
\end{para}

\begin{para}
\label{statement}
  Let $H$ be as in Proposition \ref{p:linalg}. 
  We have identified $H$ with a stalk endowed with the action of $\pi_1^{\log}(X)$. 
  Then the assumption in there means that the action of $\gamma_1$ is unipotent, and the conclusion in there is the equivalent form (iv) in Part I, Proposition 6.11 (3), which we will prove. 

  To prove this, it is enough to show $N_1=0$ on $H$.
  For, it implies $\gamma_1=\gamma_2=1$ on $H$, and hence, by \ref{pi1fiber}, 
$\pi_1^{\log}(X\times_Y \bar y(\log))$ acts trivially on $H$. 
  Therefore, by \ref{pi1}, $H$ is the pullback of an object of $\cB_Y$. 
  By the argument in \cite{KNU} 6.9 (to take a strict point of $X$), we see that the latter object belongs to $\cA_Y$. 
\end{para}

\begin{para}
  Now we have $FN_1F^{-1}=N_1+\frac ba N_2$.
  Replacing $N_1$ by $N_1-\frac b{a(q-1)}N_2$, we may assume $FN_1=N_1F$ and keep all the other relations: $[N_0,N_2]=[N_1,N_2]=0$, $FN_0=qN_0F$, 
$FN_2=qN_2F$, and $[N_0, N_1]=N_2$.
  It is enough to show that this new $N_1$ is zero.  
  In the rest, $N_1$ denotes the new one. 
\end{para}

\begin{para}
  To prove $N_1=0$, we use the commutativity of $N_2$ and $N_0$ and that of $N_2$ and $N_1$. 
  By $[N_0, N_1]=N_2$, they are respectively described as follows. 

(1) $N_0^2 N_1 - 2 N_0 N_1 N_0 + N_1 N_0^2 = 0$

(2) $N_1^2 N_0 - 2 N_1 N_0 N_1 + N_0 N_1^2 = 0$
\end{para}

\begin{para}
\label{observation}
  The proof is based on the following observation. 
  In the following, the weight means the Frobenius weight. 
  The Frobenius weight gives a splitting of the (shifted) relative monodromy filtration of $N_0$, and $N_0$ is of weight $-2$ with respect to this splitting by Part I, Lemma 4.2 (1).
  Let $H'$ be a vector subspace of $H$.
  Assume that there is a base of $H'$ of the form $\{N_0^i(h_j)\}_{i,j}$, 
where, for each $j$, $h_j$ is an $F$-eigenvector which is $N_0$-primitive and 
$N_0^i(h_j)$ are all nonzero iterated images of it. 
  Then there is a decomposition $H'=\bigoplus_w \gr_w(H')$. 
  Assume also that each $\gr_w(H')$ is closed under $N_1$. 
  Assume further $N_2=0$ on $H'$. 
  Then $H'$ is a subobject of $H$ in $\cA_X$. 

  To see it, first, since $N_1$ is of weight 0 by $FN_1=N_1F$, $H'$ is closed under the action of $N_1$. 
  Since $N_2=0$ on $H'$, it is closed under the action of $\pi_1^{\log}(X\times_Y \bar y(\log))$.
  Further, it is closed under the action of $F$ because 
all $N^i_0(h_j)$ are $F$-eigenvectors by $FN_0=qN_0F$. 
  Hence it is closed under the action of $\pi_1^{\log}(X)$.
  The admissibility and the purity 
(the conditions in Part I, 1.1) 
are satisfied at any nonsingular closed point because the relative weight filtration with respect to $N_0$ on $H'$ coincides with the restriction of that on $H$. 
  They are satisfied also at the singular point by the same reason because $N_2=0$ on $H'$. 
\end{para}

\begin{para}
We focus the Jordan normal form of $N_0$.

(3) We may assume that the sizes of the blocks are all even or all odd. 

This is because the direct sum of all subspaces corresponding to blocks of odd (resp.\ even) size is a subobject.
  In fact, it is a subobject as an $\ell$-adic sheaf by a similar argument in \ref{observation}, and, since the action of the local monodromy cone at any closed point 
on $H$ is the direct sum of the actions on the two subspaces, 
the admissibility and the purity are satisfied for each one of the two subspaces. 

(4) If there is only one block, the conclusion of Proposition \ref{p:linalg} is valid. 

In fact, then, every $\gr_w$ is of dimension one, on where $N_1$ is zero. 
  Hence $N_1=0$ on $H$. 
\end{para}

  We will apply the following linear-algebraic lemma on common eigenvectors.

\begin{lem}\label{eigenv1} Let $V$ be a  finite dimensional non-zero vector space over an algebraically closed field $F$ of characteristic $0$ and let $A, B :V\to V$ be linear maps. Let $a,b,c,d\in F$, assume $aA^2+bAB+cBA+dB^2=0$, and assume that one of {\rm (i)}, {\rm (ii)}, {\rm (iii)} below is satisfied. Then $A$ and $B$ have a common eigenvector.

{\rm (i)} $4ad=(b+c)^2$, $b\neq c$,  and at least one of $a,d$ is not $0$.

{\rm (ii)} $A$ is not invertible, $b\neq 0$, $d=0$. 

{\rm (iii)} $B$ is not invertible, $c\neq 0$, $a=0$. 

\end{lem}

\begin{pf} First, assume (i). By $4ad=(b+c)^2$, there are $\alpha, \beta \in F$ such that $\alpha^2=a$, $\beta^2=d$, $2\alpha\beta=b+c$. Since at least one of $a, d$ is non-zero, there are $\gamma, \delta\in F$ such that $b= \alpha\beta+\alpha \delta- \beta \gamma$. Then $c= \alpha\beta -\alpha\delta +\beta \gamma$. We have
$$XY-YX=Y^2, \quad \text{where}\;\; X= \gamma A+\delta B, \; Y=\alpha A+\beta B.$$
Since $\alpha \delta - \beta \gamma = (b-c)/2\neq 0$, $A$ and $B$ are linear combinations of $X$ and $Y$. Hence it is sufficient to prove that $X$ and $Y$ have a common eigenvector. By induction on $n$, we have

(i)${}'$ \, $XY^n-Y^nX=nY^{n+1}$ for all $n\geq 0$. 

Since  $\text{Tr}(XY^n-Y^nX)= 0$, we have $\text{Tr}(Y^n)=0$ for all $n>0$, and this shows that $Y$ is nilpotent. Let $v$ be an eigenvector of $X$  and take $n\geq 0$ such that $Y^nv\neq0$ and $Y^{n+1}v=0$. Then by (i)${}'$, $Y^n v$ is a common eigenvector of $X$ and $Y$. 

Next, assume (ii). Then $B$ induces $\Ker(A)\to \Ker(A)$ and hence we can find a common eigenvector  in $\Ker(A)$. Swapping $A$ and $B$ here, we obtain the proof for the case (iii).
\end{pf}

\begin{rem}\label{eigenv2}  In \ref{pflinalg}, we will use  the following cases (i) and (ii) of Lemma \ref{eigenv1}. For linear maps $A,B:V\to V$ on a finite dimensional non-zero vector space $V$ over an algebraically closed field of characteristic $0$, we have:

(i) If $A^2-2AB+B^2=0$, $A$ and $B$ have a common eigenvector. (A special case of (i) in Lemma \ref{eigenv1}.)

(ii) If $AB=-BA$ and $A$ is not invertible, $A$ and $B$ have a common eigenvector.  (A special case of (ii) and also of  (iii) in Lemma \ref{eigenv1}.)

Lemma \ref{eigenv1} treats subvarieties of the affine $4$-space $\mathbb{A}^4=\{(a,b,c,d)\}$ on which $A$, $B$ (with some extra conditions) such that $aA^2+bAB+cBA+dB^2=0$ have a common eigenvector. There may be  many other such subvarieties than those treated in Lemma \ref{eigenv1}. 
\end{rem}

\begin{para}\label{pflinalg}
  We start to prove Proposition \ref{p:linalg}.
  We may and will assume that the weight of $H$ is zero.

  For any $w \in \bZ$, let $N_{1,w}$ be the part of $N_1$ on $\gr_w$. 

  By (4), we may assume that there are at least two Jordan blocks of $N_0$. 
  By (3), the sizes of Jordan blocks of $N_0$ are all even or all odd.
  First we prove the former case. 
  In this case, we have an isomorphism 
$N_0:\gr_1\overset{\cong}\to \gr_{-1}$. 
  Let $V$ be an $F$-eigenspace of $\gr_1$. 
  Let $\lambda$ be its eigenvalue. 
  Then $N_0$ maps $V$ onto the $F$-eigenspace of $\gr_{-1}$ whose eigenvalue is $q\lambda$, which we identify with $V$.
  Then $N_{1,1}$ and $N_{1,-1}$ induce endomorphisms on $V$ because $F$ and $N_1$ commute, which we denote by the same symbols. 
  Then (2) implies 

$N_{1,-1}^2  - 2 N_{1,-1} N_{1,1} + N_{1,1}^2 = 0$ on $V$.

  By Remark \ref{eigenv2} (i), there is a common eigenvector in $V$ of $N_{1,\pm 1}$, which means that there is a nonzero vector $v \in V \subset \gr_1$ such that $N_1(v)=0$ and 
$N_1N_0(v)=0$. 
  Then by (1), $N_1N_0^k(v)=0$ for any $k \ge 0$. 
  Further, let $v_w$ be the highest weight vector whose image by some power of $N_0$ coincides with $v$. 
  Let $w$ be the weight of $v_w$: $N_0^{(w-1)/2}(v_w)=v$. 
  Note that $v_w$ is also an $F$-eigenvector because $v$ is so. 
  Then again by (1), we can show inductively $N_1N_0^k(v_w)=0$ for any $k \ge 0$ as follows. 
  We already know it for $k \ge (w-1)/2$. 
  Let $k_0 \le (w-1)/2$. 
  Assume that we know it for $k \ge k_0$.
  Then by (1), $N_0^2N_1N_0^{k_0-1}(v_w)=0$. 
  Since $N_0^2$ is injective on $\gr_{w-2(k_0-1)}$, this implies 
$N_1N_0^{k_0-1}(v_w)=0$. 
  Thus, by \ref{observation} (take $v_w$ as $h_j$ in \ref{observation}), we find a subobject generated by $N_0^k(v_w)$ for all $k \ge 0$.  
  (Note that, on this subspace, $N_1=0$ and $N_1N_0=0$ so that $N_2=[N_0,N_1]=0$, i.e., the last condition in \ref{observation} is satisfied.) 
  Since there are at least two Jordan blocks, it is a proper subobject, which is a contradiction. 

  We prove the latter case.
  In this case, we have an isomorphism $N_0^2:\gr_2 \overset\cong \to 
\gr_{-2}$.
  When $\gr_2=\{0\}$, any common eigenvector of $N_{1,0}$ and $F$ 
generates a 1-dimensional subobject. 
  In the rest, we assume $\gr_2\not=\{0\}$. 
  Let $V$ be an $F$-eigenspace of $\gr_2$. 
  Let $\lambda$ be its eigenvalue. 
  Then $N_0^2$ maps $V$ onto the $F$-eigenspace of $\gr_{-2}$ whose eigenvalue is $q^2\lambda$, which we identify with $V$.
  On the other hand, $V$ is regarded as a subspace of 
the $F$-eigenspace $\tilde V$ of $\gr_0$ whose eigenvalue is $q\lambda$ 
via $N_0:\gr_2 \to \gr_0$.  
  Further, let $p\colon \tilde V \to V$ be the section induced by $N_0:\gr_0 \to \gr_{-2}$. 
  Then $N_{1,2}$, $pN_{1,0}$, $pN_{1,0}^2$, and $N_{1,-2}$ are regarded as endomorphisms on $V$. 
  Then (1) and (2) imply 

$(1)'$ $N_{1,2} - 2p N_{1,0} + N_{1,-2} =0$.

$(2)'$ $pN_{1,0}^2  - 2 pN_{1,0} N_{1,2} + N_{1,2}^2 = 0$, and

$(2)''$ $N_{1,-2}^2  - 2 N_{1,-2} pN_{1,0} + pN_{1,0}^2 = 0$ on $V$. 

\smallskip

\noindent  Claim.  There is a common eigenvector of $N_{1,2}$ and $pN_{1,0}$. 

\smallskip

  This claim implies the proposition.
  In fact, this claim means that there is a nonzero vector $v \in V \subset \gr_2$ such that $N_1(v)=0$ and $N_0N_1N_0(v)=0$. 
  If $N_1N_0(v)$ is not zero, then, since  $N_1^2N_0(v)=0$ by (2), $N_1N_0(v)$ alone generates a 1-dimensional subobject. 
  Otherwise, $N_1N_0(v)=0$.  
  Then exactly in the same way as in the former case, we can find a proper subobject generated by $N_0^k(v_w)$ for all $k \ge0$, where 
$v_w$ is the highest weight vector whose image by some power of $N_0$ coincides with $v$. 

  We prove Claim. 
  By Remark \ref{eigenv2} (ii), it is enough to show $A:=N_{1,2}$ and $B:=pN_{1,0}$ satisfy $AB+BA=0$. 
  Let $C:=N_{1,-2}$. 
  By $(2)'$ and $(2)''$, we have $pN_{1,0}^2-2BA+A^2=0$ and 
$C^2-2CB+pN_{1,0}^2=0$. 
  Substitute $C=2B-A$ by $(1)'$ to the last equality, we have 
\begin{align*}
(2B-A)^2-2(2B-A)B+pN_{1,0}^2&=
4B^2-4BA-4AB+A^2-4B^2+2AB+pN_{1,0}^2\\
&=pN_{1,0}^2-4BA-2AB+A^2=0. 
\end{align*}
  Together with $pN_{1,0}^2-2BA+A^2=0$, we have $AB+BA=0$, as desired. 
  This completes the proof of Proposition \ref{p:linalg}.
\end{para}

\noindent {\rm Kazuya KATO
\\
Department of mathematics
\\
University of Chicago
\\
Chicago, Illinois, 60637, USA}
\\
{\tt kkato@uchicago.edu}

\bigskip

\noindent {\rm Chikara NAKAYAMA
\\
Department of Economics 
\\
Hitotsubashi University 
\\
2-1 Naka, Kunitachi, Tokyo 186-8601, Japan}
\\
{\tt c.nakayama@r.hit-u.ac.jp}

\bigskip

\noindent
{\rm Sampei USUI${}^{\sharp}$
\\
Graduate School of Science
\\
Osaka University}
\\
\footnote[0]{${\sharp}$ \ \ the deceased}
\end{document}